\documentclass{article}
\usepackage{amssymb,amsmath}
\usepackage[mathscr]{euscript}
\input{epsf}

\newcounter{sec}

\newcounter{punct}[sec]

\def\punct{\refstepcounter{punct}{\arabic{sec}.\arabic{punct}.  }}

\newtheorem{theorem}{Theorem}[sec]
\newtheorem{proposition}[theorem]{Proposition}

\newtheorem{lemma}[theorem]{Lemma}

\newtheorem{corollary}[theorem]{Corollary}

\def\COUNTERS{\addtocounter{sec}{1}
              \setcounter{punct}{0}
          \setcounter{equation}{0}
          \setcounter{theorem}{0}
          }
          
          \def\sm{\smallskip}

\begin{document}

\newcommand{\supp}{\mathop {\mathrm {supp}}\nolimits}
\newcommand{\rk}{\mathop {\mathrm {rk}}\nolimits}
\newcommand{\Aut}{\mathop {\mathrm {Aut}}\nolimits}
\newcommand{\Out}{\mathop {\mathrm {Out}}\nolimits}
\renewcommand{\Re}{\mathop {\mathrm {Re}}\nolimits}
\newcommand{\Mor}{\mathop {\mathrm {Mor}}\nolimits}
\newcommand{\Skel}{\mathop {\mathrm {Skel}}\nolimits}
\newcommand{\im}{\mathop {\mathrm {im}}\nolimits}

\newcommand{\GL}{\mathop {\mathrm {GL}}\nolimits}
\newcommand{\Sp}{\mathop {\mathrm {Sp}}\nolimits}
\newcommand{\SO}{\mathop {\mathrm {SO}}\nolimits}
\renewcommand{\O}{\mathop {\mathrm {O}}\nolimits}
\newcommand{\U}{\mathop {\mathrm {U}}\nolimits}
\newcommand{\SU}{\mathop {\mathrm {SU}}\nolimits}
\newcommand{\Ad}{\mathop {\mathrm {Ad}}\nolimits}
\newcommand{\Isom}{\mathop {\mathrm {Isom}}\nolimits}
\newcommand{\Abs}{\mathop {\mathrm {Abs}}\nolimits}
\newcommand{\Clop}{\mathop {\mathrm {Clop}}\nolimits}
\newcommand{\spike}{\mathop {\mathrm {spike}}\nolimits}
\renewcommand{\vert}{\mathop {\mathrm {vert}}\nolimits}
\newcommand{\edge}{\mathop {\mathrm {edge}}\nolimits}
\newcommand{\Hie}{\mathop {\mathrm {Hier}}\nolimits}
\newcommand{\SL}{\mathop {\mathrm {SL}}\nolimits}
\newcommand{\PSL}{\mathop {\mathrm {PSL}}\nolimits}
\newcommand{\Ba}{\mathop {\mathrm {Ba}}\nolimits}
\newcommand{\Br}{\mathop {\mathrm {Br}}\nolimits}
\newcommand{\Diff}{\mathop {\mathrm {Diff}}\nolimits}

\def\Ss{\mathrm{S}}

\def\ov{\overline}
\def\wh{\widehat}
\def\wt{\widetilde}

\renewcommand{\rk}{\mathop {\mathrm {rk}}\nolimits}
\renewcommand{\Aut}{\mathop {\mathrm {Aut}}\nolimits}
\renewcommand{\Re}{\mathop {\mathrm {Re}}\nolimits}
\renewcommand{\Im}{\mathop {\mathrm {Im}}\nolimits}
\newcommand{\sgn}{\mathop {\mathrm {sgn}}\nolimits}
\newcommand{\PGL}{\mathop {\mathrm {PGL}}\nolimits}
\newcommand{\tr}{\mathop {\mathrm {tr}}\nolimits}

\def\GLO{\mathbf{GL}}

\def\bfa{\mathbf a}
\def\bfb{\mathbf b}
\def\bfc{\mathbf c}
\def\bfd{\mathbf d}
\def\bfe{\mathbf e}
\def\bff{\mathbf f}
\def\bfg{\mathbf g}
\def\bfh{\mathbf h}
\def\bfi{\mathbf i}
\def\bfj{\mathbf j}
\def\bfk{\mathbf k}
\def\bfl{\mathbf l}
\def\bfm{\mathbf m}
\def\bfn{\mathbf n}
\def\bfo{\mathbf o}
\def\bfp{\mathbf p}
\def\bfq{\mathbf q}
\def\bfr{\mathbf r}
\def\bfs{\mathbf s}
\def\bft{\mathbf t}
\def\bfu{\mathbf u}
\def\bfv{\mathbf v}
\def\bfw{\mathbf w}
\def\bfx{\mathbf x}
\def\bfy{\mathbf y}
\def\bfz{\mathbf z}

\def\bfA{\mathbf A}
\def\bfB{\mathbf B}
\def\bfC{\mathbf C}
\def\bfD{\mathbf D}
\def\bfE{\mathbf E}
\def\bfF{\mathbf F}
\def\bfG{\mathbf G}
\def\bfH{\mathbf H}
\def\bfI{\mathbf I}
\def\bfJ{\mathbf J}
\def\bfK{\mathbf K}
\def\bfL{\mathbf L}
\def\bfM{\mathbf M}
\def\bfN{\mathbf N}
\def\bfO{\mathbf O}
\def\bfP{\mathbf P}
\def\bfQ{\mathbf Q}
\def\bfR{\mathbf R}
\def\bfS{\mathbf S}
\def\bfT{\mathbf T}
\def\bfU{\mathbf U}
\def\bfV{\mathbf V}
\def\bfW{\mathbf W}
\def\bfX{\mathbf X}
\def\bfY{\mathbf Y}
\def\bfZ{\mathbf Z}

\def\frD{\mathfrak D}
\def\frH{\mathfrak H}
\def\frK{\mathfrak K}
\def\frL{\mathfrak L}
\def\frS{\mathfrak S}
\def\frT{\mathfrak T}
\def\frU{\mathfrak U}
\def\frV{\mathfrak V}
\def\frW{\mathfrak W}
\def\frQ{\mathfrak Q}
\def\hier{\mathfrak{Hi}}

\def\frg{\mathfrak g}
\def\frz{\mathfrak z}
\def\frm{\mathfrak m}

\def\bfw{\mathbf w}

\def\R {{\mathbb R }}
 \def\C {{\mathbb C }}
  \def\Z{{\mathbb Z}}
  \def\H{{\mathbb H}}
\def\K{{\mathbb K}}
\def\N{{\mathbb N}}
\def\Q{{\mathbb Q}}
\def\A{{\mathbb A}}

\def\T{\mathbb T}
\def\P{\mathbb P}

\def\G{\mathbb G}

\def\cA{\EuScript A}
\def\cD{\EuScript D}
\def\cL{\EuScript L}
\def\cK{\EuScript K}
\def\cM{\EuScript M}
\def\cN{\EuScript N}
\def\cR{\EuScript R}
\def\cW{\EuScript W}
\def\cY{\EuScript Y}
\def\cZ{\EuScript Z}
\def\cF{\EuScript F}
\def\cT{\EuScript T}
\def\cB{\EuScript B}
\def\cE{\EuScript E}
\def\cO{\EuScript O}
\def\cP{\EuScript P}
\def\cH{\EuScript H}

\def\bbA{\mathbb A}
\def\bbB{\mathbb B}
\def\bbD{\mathbb D}
\def\bbE{\mathbb E}
\def\bbF{\mathbb F}
\def\bbG{\mathbb G}
\def\bbI{\mathbb I}
\def\bbJ{\mathbb J}
\def\bbL{\mathbb L}
\def\bbM{\mathbb M}
\def\bbN{\mathbb N}
\def\bbO{\mathbb O}
\def\bbP{\mathbb P}
\def\bbQ{\mathbb Q}
\def\bbS{\mathbb S}
\def\bbT{\mathbb T}
\def\bbU{\mathbb U}
\def\bbV{\mathbb V}
\def\bbW{\mathbb W}
\def\bbX{\mathbb X}
\def\bbY{\mathbb Y}

\def\kappa{\varkappa}
\def\epsilon{\varepsilon}
\def\phi{\varphi}
\def\le{\leqslant}
\def\ge{\geqslant}

\def\B{\mathrm B}

\def\la{\langle}
\def\ra{\rangle}

\def\F{{}_2F_1}
\def\FF{{}^{\vphantom{\C}}_2F_1^\C}

\newcommand{\dd}[1]{\,d\,{\overline{\overline{#1}}} }

\def\lambdA{{\boldsymbol{\lambda}}}
\def\alphA{{\boldsymbol{\alpha}}}
\def\betA{{\boldsymbol{\beta}}}
\def\mU{{\boldsymbol{\mu}}}
\def\PI{{\boldsymbol{\Pi}}}

\def\1{\boldsymbol{1}}
\def\2{\boldsymbol{2}}

\def\Th{\mathrm{T\!h}}
\def\T{\bbT}
\def\B{\bbB}

\def\ls{(\!(}
\def\rs{)\!)}
\def\Ls{\bigl(\!\bigr(}
\def\Rs{\bigl)\!\bigr)}
\def\LS{\Bigl(\!\!\Bigr(}
\def\RS{\Bigl)\!\!\Bigr)}

\def\tochka{\boldsymbol{\cdot}}

\def\black{\mathrm{black}}
\def\blue{\mathrm{blue}}
\def\red{\mathrm{red}}
\def\fin{\mathrm{fin}}

\def\procent{}

\begin{center}
\bf\Large

On the group of spheromorphisms of
\\ the homogeneous non-locally finite tree

\bigskip

\large \sc Yury A. Neretin%
\footnote{The research was supported by the grants FWF, Projects  P28421, P31591.}
\end{center}

{\small 
	Consider a tree $\T$, all whose vertices have countable valence;
	its boundary is the Baire space $\B\simeq\N^\N$; 
	continued fractions expansions identify the set of irrational numbers
	$\R\setminus \Q$  with $\B$.
	Removing $k$ edges from $\T$
	we get a forest  consisting of copies of $\T$. 
	A spheromorphism (or hierarchomorphism)
	of $\T$ is an isomorphisms of two
	such subforests considered as a transformation of $\T$
	or of $\B$. Denote the group of all spheromorphisms by $\Hie(\T)$.
	We a show that the   correspondence $\R\setminus\Q\simeq\B$
	sends
	the Thompson
	group
	realized by piecewise $\PSL_2(\Z)$-transformations   to a subgroup of $\Hie(\T)$.  
	We construct some unitary representations of the group
	$\Hie(\T)$, show that the group
	of automorphisms $\Aut(\T)$ is spherical in $\Hie(\T)$,
	and describe the train (enveloping category)
	of $\Hie(\T)$.
}

\section{Introduction}

\COUNTERS

{\bf \punct The tree $\T$ and its boundary $\partial\T$.%
\label{ss:1.1}}
For a set $A$ denote by $\# A$ the number of its elements.
Denote by $\Z_+$ the set of all nonnegative integers.

 Recall that a {\it tree} $T$ is a connected graph without cycles.
 Denote by $\vert(T)$ the set of its vertices, $\edge(T)$ the set 
 of its edges. A {\it forest} is a disjoint union of trees.
  We admit both finite and infinite trees.

Denote by $\T$ the tree such that each vertex is contained in a countable number
of edges. Such a tree is unique up to isomorphisms of trees. It can be realized
in the following form. Vertices of $\T$ are enumerated
by finite collections
\begin{equation}
(s_1, s_2,\dots,s_m),\qquad \text{where $s_1\in \Z_+$,
	$s_2$, \dots, $s_m\in \N$,}  
\label{eq:ssss}
\end{equation}
where $m=0$, 1, 2, \dots. We admit an empty collection,
 below we call such vertex the {\it initial point%
	\footnote{We do not use the term 'root' since we regard $\T$ as a non-rooted
		tree.}
	 of $\T$}
and denote by '$\tochka$'. Edges have form
$$
(s_1,\dots, s_m)\, \text{---}\, (s_1,\dots, s_m,s_{m+1}).
$$

 We say that a {\it way} in $\T$ is a sequence of pairwise distinct vertices 
$a_0$, $a_1$, $a_2$, \dots such that $a_j$ and $a_{j+1}$ are connected by an edge.
We say that two ways $a=\{a_0, a_1,a_2,\dots \}$,
$b=\{b_0, b_1,b_2,\dots \}$  are {\it equivalent} if there is $k\in\Z$
such that for sufficiently large $j$ we have $b_j=a_{j+k}$.
The {\it boundary} $\partial\T$ of $\T$ is the set of all ways defined up to this equivalence.
	Fix a vertex $r$. Then for any point $\omega\in \partial \T$
there is a unique way starting at $r$ and coming to $\omega$
(formally, the last phrase means that there is a unique representative 
of $\omega$ starting at $r$).	
 Define a distance
	between ways $a=\{r, a_1,a_2,\dots \}$ 
	$b=\{r, b_1,b_2,\dots \}$ by
	$\mathrm{dist}_r(a,b)=2^{-j}$, where $j$ is the first number 
	 such that $a_j\ne b_j$. Then $\partial\T$
	becomes a complete totally disconnected metric space. 
	Distances $\mathrm{dist}_r(\cdot,\cdot)$ depend on $r$ but they
	define the same topology.
	
	Choosing $r=\tochka$,
	 we identify $\partial\T$ with the set of all sequences
	\begin{equation*}
	(s_1, s_2, s_3\dots),\qquad \text{where $s_1\in \Z_+$,
		$s_2$, $s_3, \dots\in \N$,}  
	\end{equation*}

\sm

{\bf\punct The Baire space and continued fractions.}
Recall (see, e.~g., \cite{Kech}) that the {\it Baire space}  $\B$ is the topological space
homeomorphic to the countable product of countable discrete spaces,
$$
\B\simeq \N^\N=\N\times \N\times \N\times \dots
$$
equipped with the Tikhonov topology. Clearly,
 the boundary $\partial\T$ is
 $\B$.

The Baire space can be identified with the set $\R\setminus \Q$
of irrational numbers. Namely, let $x\in\R\setminus \Q$.
Consider its
 continued fraction decomposition,
$$
x=u_0+\cfrac{1}{u_1+\cfrac{1}{u_2+\cfrac{1}{u_3+\cfrac{1}{u_4+\dots}}}}
=:[u_0; u_1,u_2,u_3,\dots].
$$
For irrational $x$ the continued
 fraction is infinite, therefore we have an identification
$$
\R\setminus\Q\simeq \Z\times \N\times \N\times \N\times\dots\simeq\B
$$
The Baire space $\B$ and this correspondence had a fundamental role
in works of Luzin on descriptive set theory in 1920s, see \cite{Lusin0}, \cite{Lusin}.

\sm 

{\bf \punct The group of spheromorphisms.%
\label{ss:spheromorphisms}}
Denote by $\Aut(\T)$ the {\it group of all automorphisms
	of $\T$}. We define the topology on $\Aut(\T)$ assuming that all point-wise
 stabilizers
$\cK(J)$
 of finite subtrees $J\subset \T$ are open subgroups in $\Aut(\T)$.
We get a Polish group%
\footnote{i.~e., a topological group that is homeomorphic
	(as a topological space) to a complete metric space, see, e.g.?
	\cite{Kech}.}.

\begin{figure}
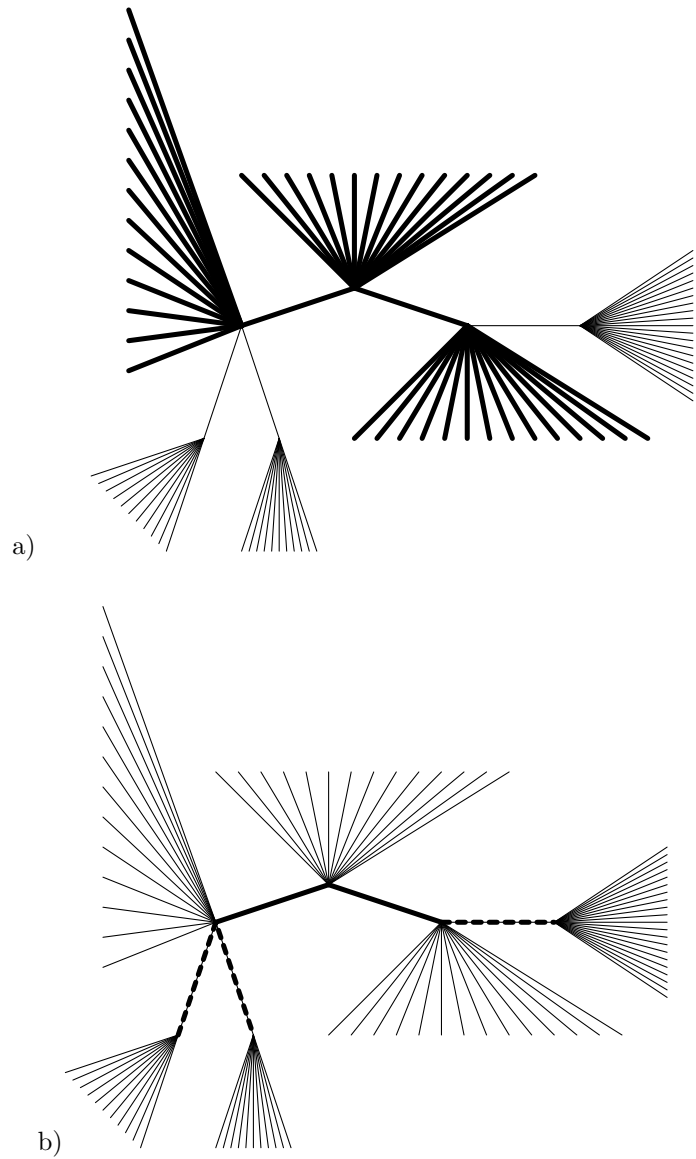

$$\mathrm{a)}\qquad \epsfbox{inftree.1}  $$
	
$$\mathrm{b)} \epsfbox{inftree.2}  $$
	
\caption{
	Refs. to Subsect. \ref{ss:spheromorphisms}
	and \ref{ss:frame}.
	\newline
	a) A $(\T)$-subtree.
\newline
b) The corresponding frame.}
\label{fig:T-subtree}	
\end{figure}

 Consider a proper subtree $S\subset \T$
isomorphic to $\T$. We say that $S$ is a {\it $(\T)$-subtree}
if there is only finite number of edges $[a,b]$ of
$\T$ such that $a\in S$ and $b\notin S$, see Fig. \ref{fig:T-subtree}.a.
An intersection of two $(\T)$-subtrees is $(\T)$-subtree
or the empty set.
If $(\T)$-subtrees $P$, $Q$ have a common vertex, then
$P\cup Q$ is a $(\T)$-subtree. 

We say that a $(\T)$-{\it covering forest} of $\T$ is a finite collection
of disjoint $(\T)$-subtrees $S_1$, \dots, $S_k$
such that 
$$
\vert(\T)=\bigcup \vert(S_j),\qquad
 \text{the set $\edge(\T)\setminus \bigcup \edge(S_j)$ 
 	is finite.}
$$
In other words, a $(\T)$-covering forest is obtained from $\T$ by a removing
 a finite collection of edges. 

\sm 

Let $S_1$, \dots, $S_k$ and $R_1$, \dots, $R_k$ be two $(\T)$-covering
forests. A {\it spheromorphism} (or {\it hierarchomorphism})
$g$
of $\T$ is
a collection $\{g^{(j)}\}$ of isomorphisms
$g^{(j)}:S_j\to R_j$.
 Notice that a spheromorphism $g$ determines 
 a permutation
 \begin{equation}\vert(\T)\to \vert(\T)
 \label{eq:vert-vert}
 \end{equation}
 and a homeomorphism
\begin{equation} \partial\T\to\partial\T.
 \label{eq:T-T}
\end{equation}
 These maps determine one another, 
 two spheromorphisms are equal if the
 coresponding  maps (\ref{eq:vert-vert}) (or, equivalently, (\ref{eq:T-T}))
 coincide,  we denote them by the same symbol $g$.

Let $\{g^{(j)}\}:\{S^{(j)}\}\to \{R^{(j)}\}$,
$\{h^{(l)}\}:\{Q^{(l)}\}\to \{T^{(l)}\}$
be spheromorphisms. 
  Their product 
 is the collection of isomorphisms of $(\T)$-subtrees.
 $$
 \{h^{(l)} g^{(j)}\}:
 \{(g^{(j)})^{-1} (R^{(j)} \cap Q^{(l)}) \}
 \to \{h^{(l)} (R^{(j)} \cap Q^{(l)}) \}.
 $$
 Denote by $\Hie(\T)$ the group of all spheromorphisms.
 By definition, the group $\Hie(\T)$
 is embedded to the group of all permutations
 of $\vert(\T)$ and the group of all homeomorphisms
 of $\partial\T$. 
 
\sm


\sm 

We define a {\it topology} on $\Hie(\T)$ from two conditions:

\sm

a) the induced topology on $\Aut(\T)\subset \Hie(\T)$
coincides with the natural topology on $\Aut(\T)$.

\sm 

b) The topology on the countable homogeneous space 
$\Hie(\T)/\Aut(\T)$ is discrete.

\sm 

In this way we get a structure of a Polish group on $\Hie(\T)$.

\sm


\sm

{\bf \punct Thompson group and the group of spheromorphisms.%
\label{ss:thompson}}
Consider the natural action of the group $\PGL_2(\Z)$
on the real projective line $\R\P^1$ by linear fractional transformations, 
$$
x\mapsto \frac{ax+b}{cx+d}, \qquad \text{where $a$, $b$, $c$, $d\in \Z$,
	$ad-bc=\pm1$.} 
$$
It contains a subgroup $\PSL_2(\Z)$ consisting of transformations
with $ad-bc=1$.
Clearly, the set $\Q$ of rational numbers is invariant
with respect to such transformations,
 therefore $\PGL_2(\Z)$ acts also on the set of irrational numbers $\R\setminus \Q$.
Next, consider the {\it Thompson group} $\Th$ of all continuous piece-wise $\PSL_2(\Z)$-transformations.
It is easy to show that such transformations have smoothness $C^1$
and points of break of the second derivative are rational. 

More constructive description of such transformations 
is given on Fig. \ref{fig:modular}.

\begin{figure}
		\epsfbox{haar.2}, \qquad \epsfbox{haar.5}
	\caption{3
		Refs. to Subsect. \ref{ss:thompson}.
		\newline
		a) Consider the Lobachevsky plane $\Lambda:\Im z>0$. 
		Lines on $\Lambda$ are semicircles and rays orthogonal to the line $\Im z=0$.
		Consider the line
	 $\Re z=0$ and all its images under elements of the group $\PSL_2(\Z)$.
	We get a family $\cZ$ of lines, each line separates $\Lambda$
	into two half-planes (on the Figure they are half-disks, 
	complements to half-disks, or right angles).
	For any pair of two such half-planes there is a unique element
	of $\PSL_2(\Z)$ sending one half-plane to another.
		\newline
		b) 
		 We take two ideal $n$-gons $U$, $V$, whose sides are 
		 contained in the family $\cZ$.  Let $A_1$, \dots, $A_n$ be complementary
		half-planes  to
		$U$  and $B_1$, \dots, $B_n$
 to $V$ (enumerated in according to the natural cyclic order).
		Consider a cyclic permutation
		$j\mapsto k+j (\mod n)$ of the set
		$\{1,2,\dots,n\}$.
		For each $j$ consider the canonical $\PSL(2,\Z)$-transformation
		$A_j\mapsto B_{j+k}$. In this way we get a piece-wise
		$\PSL_2(\Z)$-transformation of $\R\P^1$. 
	}
	\label{fig:modular}
\end{figure}

\sm

{\sc Remark.} The Thompson group $\Th$ was defined by Richard Thompson
in 1966 as a counterexample, later it became clear that it is
a very interesting discrete group with unusual properties,
see, e.~g., \cite{GS}, \cite{Imb}, \cite{Fos}. May be the most
strange is its relation with the Minkowski function $?(x)$
(apparently observed by Sergiescu, see \cite{Fos}).
\hfill $\boxtimes$

\sm

The group $\Th$ acts on $\R\setminus \Q$ and therefore acts on the Baire space
$\partial\T$. Theorem \ref{th:thompson} shows that {\it the Thompson group acts
on $\partial\T$ by spheromorphisms%
\footnote{Let $[a_0;a_1,a_2,\dots]$ be a continued fraction,
	 $[b_0;b_1,b_2,\dots]$ be its image under 
	$\PSL_2(\Z)$-transformation. Clearly, there is $m\in\Z$
	such that for sufficiently large 
	$j$ we have $b_{j+m}=a_j$. Therefore the same holds 
	for transformations from the Thompson group.
	Notice that such tail equivalence does not imply our statement.}.}

\sm

{\bf \punct The train of $\Hie(\T)$.}
%
In Section 4 we get a {\it train construction}
in the sense of \cite{Ner-book}.
Recall  that quite often a pair $G\supset K$ 
of infinite-dimensional groups generates a natural category
({\it train of  $(G,K)$}) acting in unitary representations
of $G$, such groups are called {\it $(G,K)$-pairs}.
 In our case $G=\Hie(\T)$, $K$ is the stabilizer
$\cK\subset \Aut(\T)$ of a vertex.

 Namely, let $J$
be a {\it nonempty}%
\footnote{The construction below
does not hold for double cosets
$\Aut(\T)\setminus\Hie(\T)/\Aut(\T)$.
} finite subtree in $\T$, let $\cK(J)\subset \Aut(\T)$
be its point-wise stabilizer.
In Subsect. \ref{ss:bi-tree} we give a combinatorial description of double coset spaces%
\footnote{Let $G$ be a group, $K_1$, $K_2$
subgroups. A {\it double coset} is a subset
in $G$
of the form $K_1gK_2$, where $g\in G$.
The set $K_1\setminus G/K_2$ denotes
the space of all double cosets.
\newline 
 If a subgroup
$K$ is compact, then there is a natural structure
of a 'hypergroup' on $K\setminus G/K$, i.~e.,
we have a map from $K\setminus G/K\times K\setminus G/K$
to the space $\cM(K\setminus G/K)$ of measures on
$K\setminus G/K$.  Namely, we consider uniform probabilistic
measures $\mu_{g_1}$, $\mu_{g_2}$ on double cosets $Kg_1K$, $Kg_2K$,
  decompose their convolution
 $\mu_{g_1}*\mu_{g_2}=\int_{K\setminus G/K} \mu_g \,d\nu(g)$,
 and get a probabilistic measure $\nu$ on $K\setminus G/K$ depending on $Kg_1K$,
 $Kg_2K$.
 \newline
 The situation discussed below has not analogs
 for locally compact groups and  is relatively usual for infinite dimensional
 groups, namely
 we have associative multiplications
 on sets $K\setminus G/K$, where $G=\Hie(\T)$,
 $K=\cK(J)$.}
$
\cK(J_1)\setminus \Hie(\T)/ \cK(J_2)
$ in terms of colored graphs.
 Next, we show that there is a natural
$\odot$-multiplication 
$$
\cK(J_1)\setminus \Hie(\T)/ \cK(J_2)
\,\,\times \cK(J_2)\setminus \Hie(\T)/ \cK(J_3)
\,\,\longrightarrow \,\, \cK(J_1)\setminus \Hie(\T)/ \cK(J_3),
$$
and it determines a structure of a category.
Objects of the category 
$\hier$ are nonempty finite subtrees $J\subset \T$, morphisms
are double cosets,
$$
\Mor_\hier(J,I):=\cK(I)\setminus \Hie(\T)/ \cK(J).
$$

Consider a unitary representation $\rho$ of the group $\Hie(\T)$
in a Hilbert space $V$. Denote by $V(J)\subset V$
the subspace of all $\cK(J)$-fixed vectors in $V$. Denote by
$P(J)$ the operator of orthogonal projection to $V(J)$.
For any subtrees $I$, $J$ and $g\in \Hie(\T)$ we define operators 
$$
\wt\rho_{I,J}(g)= P(I)\,\rho(g) \, P(J):\, V(J)\to V(I).
$$
This operator depends only on the double coset
$\frg=\cK(I)\, g \, \cK(J)$ containing $g$.

We show (Theorem \ref{th:multi}) that
for  $\frg_1\in \Mor(J_2, J_1)$,  $\frg_2\in \Mor(J_3, J_2)$
we have
$$
\wt \rho_{J_1,J_2}(\frg_1)\,\wt\rho_{J_2,J_3}(\frg_2)=
\wt\rho_{J_1,J_3}(\frg_1\odot\frg_2).
$$

Notice, that this operation is a representative of a huge zoo
of train constructions for $(G,K)$-pairs
(see, e.~g., \cite{Olsh-symm}, \cite{Olsh-GB}, 
\cite{Ner-book}, \cite{Ner-symm}).

\sm 

{\bf \punct Unitary representations of $\Hie(\T)$.}
We show that $\Hie(\T)\supset \Aut(\T)$
is a spherical pair, i.e., any irreducible
unitary representation of $\Hie(\T)$ has at most
one (up to a scalar factor) non-zero $\Aut(\T)$-fixed 
vector (Theorem \ref{th:sphericity}).

We also show that $\Hie(\T)$ has non-trivial $\Aut(\T)$-spherical representations.
 In fact,
we construct embeddings of $\Hie(\T)$ to two Olshanski's spherical
$(G,K)$-pairs,
the first one is related to infinite-symmetric groups
(Subsect. \ref{ss:topology}), the second is related
to classical groups (Sect.\ref{s:GL}). In both cases 
the subgroup $\Aut(\T)$ embeds to $K$,
 restricting $K$-spherical representations%
\footnote{In the both cases classification of $K$-spherical
representations of $G$ is known, see \cite{Olsh-symm}, \cite{Pick}.}
of $G$ to $\Hie(\T)$ we get $\Aut(\T)$-spherical representations
of $\Hie(\T)$. 

\sm 

{\bf \punct Similar groups.}
Consider the Bruhat--Tits tree $\cT_p$,
i.~e., the tree whose vertices have valence $p+1$.
According Bruhat and Tits, such trees are 
$p$-adic couneterparts of the Lobachevsky plane,
and more generally of noncompact rank 1 Riemannian
symmetric spaces, see e.g., \cite{Ner-gauss}, Sect. 10.4.
Applying the same approach to a non-Archimedean 
field with discrete absolute value and countable residue field%
\footnote{For instance, we can consider
	the field of formal Laurent series over
	$\Q$, then the group $\PGL_2$ over this field acts
	in a natural way by automorphisms of the tree  $\T$.}, we get the tree  $\T$.
The group $\Aut(\cT_p)$ of all automorphisms of $\cT_p$
is counterpart of real and $p$-adic groups
$\SL_2$ on the level of representation theory,
see \cite{Olsh-trees}. A {\it spheromorphism}
of $\cT_p$ is a homeomorphism $q$ of the boundary
$\partial\cT_p$ such that for any point of the boundary
there is a neighborhood, where $q$ coincides with an
automorphism of the tree. The group
$\Hie(\cT_p)$ of all spheromorphisms 
was defined in \cite{Ner1}--\cite{Ner2} as a 
counterpart of the group $\Diff(S^1)$ of diffeomorphisms
of the circle and  the group%
\footnote{In particular, there is a natural
	inclusion $\Diff(\P\Q_p^1)\subset \Hie(\cT_p)$.}
$\Diff(\P\Q_p^1)$
of locally analytic diffeomorphisms of the $p$-adic
projective line $\P\Q_p^1$.
The group 
$\Hie(\cT_p)$
has  
numerous  properties unusual for locally compact groups, see a long list of references
in \cite{Ner-19}.  

So we have a  family
of groups including 
\begin{equation}
\text{
$\Diff(S^1)$,  $\Diff(\P\Q_p^1)$,
$\Hie(\cT_p)$, $\Hie(\T)$.}
\label{eq:groups}
\end{equation}
The group $\Hie(\T)$ looks  like a monster,
however as an object of representation theory
it is simpler than its relatives. The reason
is a presence of the subgroup
$\cK(\tochka)$, which is 'heavy' in the sense of
\cite{Ner-book}. So understanding of the  group
$\Hie(\T)$
can be useful as a standpoint for investigation 
of other groups (\ref{eq:groups}).

\sm 

{\bf Acknowledgments.} The author is grateful to
Vlad Sergiescu and Nikita Gorbatyuk for discussions of this topic.

 \section{Preliminary remarks%
 	\label{s:generalities}}
 
 \COUNTERS
 
 {\bf \punct Frames of $(\T)$-subtrees.%
 \label{ss:frame}}
  For a $(\T)$-subtree $S$
 we mark all edges connecting $S$ with $\T\setminus S$. Call the {\it frame of $S$} the minimal subtree
 in $\T$ containing marked edges.
 
  Notice that terminal edges of the frame
 are precisely marked edges. If a frame has more than two
 vertices then it  uniquely determines a $(\T)$-tree. We 
 remove its terminal edges of a frame from $\T$, then $\T$ splits into a disjoint
 union of $(\T)$-trees, we choose a piece that contains
  non-terminal vertices of the frame%
  \footnote{If the frame has only two vertices,
  	then $(\T)$ splits into two parts, and we can not
  	distinguish them.}.

 Clearly, any finite subtree with $\ge 2$
 vertices can be a frame and 
  the isomorphism class of
 a frame is a unique invariant
 of $(\T)$-subtrees under the action of $\Aut(\T)$. 
 
 \sm
 
 {\sc Remark.} In \cite{Ner3} there was defined a smaller group
 of spheromorphisms $\Hie^\circ(\T)\subset \Hie(\T)$.
  Namely, we consider a $(\T)$-covering
 forest $\{S_1, \dots, S_k\}$ and a collection $g_j\in \Aut(\T)$
 such that $g_j S_j$ is a $(\T)$-covering forest.
 Then we have a spheromorphism in the sense of the previous definition.
 However, our definition allows isomorphisms $S_j\to R_j$, which have not
 extensions to automorphisms of the whole tree $\T$.
 \hfill $\boxtimes$

\sm

{\bf\punct The perfect compatible $(\T)$-forest for a spheromorphism.%
\label{ss:canonical-forest}}
We say that a {\it $(\T)$-subtree $S$ is compatible} with a spheromorphism
$g$ if the map $g:\vert(S)\to\vert(\T)$ is an embedding of trees.
We say that a {\it $(\T)$-covering forest $\{S_j\}$
is compatible with $g$} if all trees $S_j$ are compatible with
$g$.

\begin{lemma}
For any spheromorphism $g$ there is a unique compatible $(\T)$-covering forest
$P_1$, \dots, $P_l$  such that any compatible $(\T)$-covering forest
is obtained by removing  a finite collection of edges from  
the forest $\{P_i\}$.	
\end{lemma}

Let us call such forest the {\it perfect} $(\T)$-covering forest for $g$.

\sm

{\sc Proof.} Consider a compatible $(\T)$-covering forest $\{S_j\}$ with minimal possible number 
of trees, say $l$.  Let $\{Q_\alpha\}$ be another covering.
Let some $Q_\beta$ be not contained in any $S_i$. 
Consider trees 
 $S_{\gamma_1}$, $S_{\gamma_2}$, \dots, $S_{\gamma_m}$ that have nonempty intersections 
with $Q_\beta$, by definition $m\ge 2$.
Then 
$$\wt S:=Q_\beta \cup \bigl(\cup S_{\gamma_i}\bigr)
=\cup S_{\gamma_i}
$$ 
is a $(\T)$-subtree
compatible with $g$. We get a compatible $(\T)$-covering  forest with
$l-m+1<l$ elements.
\hfill $\square$

\sm

{\bf \punct Skeletons of  $(\T)$-covering forests.%
\label{ss:skeletons}}
Consider a  $(\T)$-covering forest $\{S_j\}$. Paint blue  all
edges that do not contained in the trees $S_j$. 
Consider the minimal subtree $\Sigma\subset \T$ containing
all blue edges, paint remaining edges  black. 
We call the colored tree obtained in this way  the {\it skeleton
$\Skel\{S_j\}$ of the forest}%
\footnote{The skeleton $\Skel\{S_j\}$ is union of 
frames of trees $S_j$. For frames  colorings are not necessary
since blue edges are precisely terminal edges.}
$\{S_j\}$, see Fig.  \ref{fig:skeleton}.

Clearly, any terminal vertex  of a skeleton
is contained in a blue edge. This property characterizes 
trees that can be skeletons. Orbits of $\Aut(\T)$
on the set of all $(\T)$-covering forests are enumerated
by skeletons defined up to an isomorphism.

\begin{figure}
	$$\epsfbox{inftree.3}$$
	\caption{}
	Ref. to Subsect. \ref{ss:skeletons}.
	A skeleton of a $(\T)$-covering forest.
	Blue edges are denoted by 
	$\epsfbox{inftree.4}$.
	\label{fig:skeleton}
\end{figure}

\sm

{\bf \punct Bi-trees and spheromorphisms.%
\label{ss:bi-tree1}}
Consider a finite  graph $\Gamma$
whose edges are colored black, blue, and red.
We say that $\Gamma$ is a {\it bi-tree}%
\footnote{See parallel combinatorial structures 
for the Thompson group and the groups of spheromorphisms
of Bruhat--Tits trees in \cite{Bur}, \cite{Ner-19}.}
if

\sm

$\bullet$  the subgraph $\Gamma^\black_\blue$ 
of $\Gamma$ consisting
of black and blue edges is a tree
and the same holds for the subgraph
$\Gamma^\black_\red$ 
consisting of black and red edges;

\sm

$\bullet$ $\Gamma$ has not vertices of valence 1.

\sm

In particular, the number of blue edges equals the number
of red edges; the subgraph $\Gamma^\black$
consisting of black edges is a forest.

Two bi-trees $\Gamma_1$ and $\Gamma_2$ are {\it equivalent}
if there is a color-preserving isomorphisms 
$\Gamma_1\to\Gamma_2$ of the graphs.

We wish to construct a canonical correspondence
$$
\Bigl\{\text{Set of all bi-trees} \Bigr\}
\,\longleftrightarrow\, \Aut(\T)\setminus\Hie(\T)/\Aut(\T).
$$

\sm

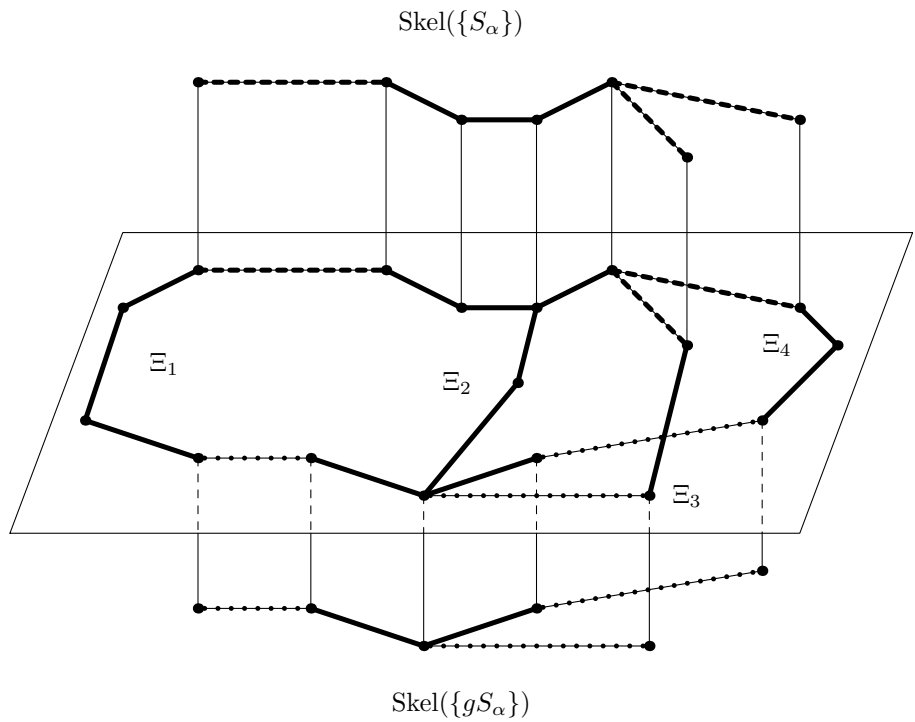
 \begin{figure}
 	$$\epsfbox{inftree.8}$$
 Notation:
 $\epsfbox{segment.6}$
  -- black,
 $\epsfbox{segment.7}$
  -- blue,
 $\epsfbox{segment.8}$
  -- red.
 Thin lines 
$\epsfbox{segment.9}$
 are auxiliary (and are not elements of graphs).
 	\caption{Ref to Subsect. \ref{ss:bi-tree1}. Skeletons $\Skel\{S_\alpha\}$
 	and $\Skel\{g S_\alpha\}$ and the corresponding
 bi-tree drawn in the 'horizontal' plane.}
\label{fig:two-skeletons}
 \end{figure}

{\sc Bi-trees of spheromorphisms}. See Fig. \ref{fig:two-skeletons}.
 Let $\{S_\alpha\}$ be the perfect
$(\T)$-covering forest for a spheromorphism $g$. 
Consider the skeleton $\Skel\{S_\alpha\}$
of $\{S_\alpha\}$, it is a tree  with black and blue edges.
Consider also the skeleton $\Skel\{gS_\alpha\}$ of $\{gS_\alpha\}$,
let us color it black and red (instead of blue). 
For each $S_j$ consider 
the minimal subtree $\Xi_j$ in $S_j$ containing the subtrees
$$
S_j\cap \Skel\{S_\alpha\},
\qquad  g^{-1}\bigl(gS_j\cap\Skel\{gS_\alpha\}\bigr)
.$$
We have embeddings of subtrees
$$
\begin{array}{ccl}
&& \Xi_j\\
&\nearrow&\\
\Xi_j\cap \Skel\{S_\alpha\}&&\\
&\searrow&\\
&&\Skel\{S_\alpha\}
\end{array},
\qquad
\begin{array}{ccl}
&& \Xi_j\\
&\nearrow&\\
\Xi_j\cap g^{-1}\Skel\{gS_\alpha\}&&\\
&\searrow&\\
&&\Skel\{gS_\alpha\}
\end{array}.
$$
We glue together trees $\Xi_1$, $\Xi_2$,
\dots, $\Skel\{S_\alpha\}$, $\Skel\{gS_\alpha\}$
identifying images of   subtrees
$\Xi_j\cap \Skel\{S_\alpha\}$ and  $\Xi_j\cap g^{-1}\Skel\{gS_\alpha\}$
 in the target-spaces
 and get a graph $\Gamma(g)$, we call it the {\it bi-tree of the spheromorphism $g$.}
 
 Formulate the last step in a simpler way.
 Consider the forest $\{\Xi_\alpha\}$.
 Consider a blue edge in $\Skel\{S_\alpha\}$. It has two ends,
 which are vertices of different $\Xi_i$, $\Xi_j$.
 We connect these vertices by a blue edge. Repeat the same
 procedure for red edges%
 \footnote{We can not receive a double blue-red edges, otherwise
 a $(\T)$-covering forest $\{S_\alpha \}$
is not perfect.}. 

\sm

{\sc Remark.} By construction, both skeletons
$\Skel\{S_\alpha \}$, $\Skel\{gS_\alpha \}$
are embedded to $\Gamma(g)$.
\hfill $\boxtimes$

\sm 
 
 {\sc Construction of a double coset
 $\in\Aut(\T)\setminus\Hie(\T)/\Aut(\T)$
 from a bi-tree $\Gamma$}.
Take two copies $\T_1$, $\T_2$
of tree $\T$. Choose isomorphisms $\theta_{1,2}:\T\to \T_{1,2}$.
Consider  embeddings 
$$p: \Gamma^\black_\blue\to \T_1,\qquad q: \Gamma^\black_\red\to \T_2.$$
Then $p$ (resp. $q$) determines a $(\T)$-covering forest
(we remove images of blue edges from $\T$), say $\{P_\alpha \}$
(say $\{Q_\alpha \}$).
The set of components of $\{P_\alpha \}$ (resp. $\{Q_\alpha \}$)  is in
a one-to-one correspondence with the set
of components $\{\Xi_\alpha\}$
 of $\Gamma^\black$. For each 
 $P_j$ we restrict $p$ to $\Gamma^\black_\blue\cap \Xi_j$
 and get an embedding of this tree to $P_j$. Extend it
 to an embedding $\wt p_j:\Xi_j\to P_j$. 
 In the same way we get embeddings $\wt q_j:\Xi_j\to Q_j$.
 Finally, we choose isomorphisms $\wt r_j:P_j\to Q_j$
 such that 
 $\wt q_j=\wt r_j\circ \wt p_j$.
 Thus we get a spheromorphism
 $r:=\{\wt r_j\}:\T_1\to \T_2$
 and define a spheromorphism
 $g:\T\to\T$ as
 $$
 g=\theta_2^{-1}\circ r\circ \theta_1.
 $$

{\bf \punct A  $(G,K)$-pair related to symmetric groups.%
\label{ss:symmetric}}
Let $\Omega$ be a countable set.
Denote by $S(\Omega)$ the group of all permutations
of $\Omega$. 
The topology on $S(\Omega)$ is determined from 
the condition: stabilizers of finite subsets in $\Omega$ are
open subgroups. This determines a structure of a Polish group
on $S(\Omega)$. On the other hand (see \cite{KR}) it is a unique
separable topology on the full infinite symmetric group
(in particular, all unitary representations of this  group
are automatically continuous in our topology).
We also write $S_\infty$ if we do not wish to indicate
the set $\Omega$.

Denote by $S^\fin(\Omega)=S^\fin_\infty$
the subgroup of finitely supported permutations of $\Omega$,
it is a countable group equipped with the discrete topology.

Let $A$ and $B$ be disjoint countable sets. Denote by
$S(A\bigl|B)$ the subgroup in $S(A\sqcup B)$ 
generated by $S(A)\times S(B)$ and $S^\fin(A\sqcup B)$.
In notation of \cite{Olsh-symm}, \cite{Ner-book}
it is a $(G,K)$-pair
$$
\bigl(S_{2\infty},S_\infty\times S_\infty).
$$
Unitary representations of this $(G,K)$-pair were classified by
Olshanski \cite{Olsh-symm}.
 
 For any element $\sigma$ of $S(A\bigl|B)$  there is a number $k$ such that
 $\sigma$ sends precisely $k$ elements of $A$ to $B$
 and $k$ elements of $B$ to $A$ (this property can be
 regarded as a definition of our group).
 

 The homogeneous space
$$\Omega:=S(A\bigl|B)\bigr/( S(A)\times S(B))$$
 is countable. It can be identified
with the set of all subsets $U\subset A\sqcup B$ such that
the sets $A\setminus U$ and $U\setminus A$ are finite
and contain the same number of elements.

We define
the topology on $S(A\bigl|B)$ from the assumptions: 

\sm

$\bullet$ the induced topology on $S(A)\times S(B)$
is the natural topology on this subgroup.

\sm 

$\bullet$ the topology on the homogeneous
space $S(A\bigl| B)\bigr/ (S(A)\times S(B))$ is discrete.

\sm 

It is clear that  we get a Polish group. 

\sm 

{\sc Remark.} The group $S(A\bigl| B)$ acts in $\ell^2(\Omega)$
and the topology of $S(A\bigl|B)$ is induced from the unitary group
of $\ell^2(\Omega)$.
\hfill $\boxtimes$



\sm

{\bf\punct The topology on $\Hie(\T)$.%
\label{ss:topology}}
The group $\Aut(\T)$ acts on the set $\vert(\T)$ of vertices
of the tree, therefore we get an embedding 
$\Aut(\T)\to S\bigl(\vert(\T)\bigr)$. The topology
on $\Aut(\T)$ defined above is induced from the symmetric group
$S\bigl(\vert(\T)\bigr)$.

Next, we define a topology on $\Hie(\T)$ from the following two conditions:

\sm 

A) the induced topology on $\Aut(\T)$ coincides
with the natural topology on  $\Aut(\T)$.

\sm 

B) this topology is a strongest topology satisfying the property
A.

\sm 

In particular, a homomorphism from $\Hie(\T)$ to a topological group
 is continuous
if and only if it is continuous on the subgroup $\Aut(\T)$.

\begin{proposition}
	{\rm a)} The topology satisfying the properties
	{\rm A--B} exists and the {\rm(}countable{\rm)} homogeneous
	space $\Hie(\T)/\Aut(\T)$ has the discrete topology.
	
	\sm 
	
	{\rm b)} The group $\Hie(\T)$ is Polish.
\end{proposition}

  To observe this, we consider the set $C$ of all non-ordered pairs 
  $(u,v)$, where $u$, $v\in\vert(\T)$, $u\ne v$.
  We put $(u,v)$ to a set $A$ if $u$, $v$ are connected by an edge
  and to $B$ otherwise. An spheromorphism $g\in\Hie(\T)$ acts
on $C$ sending $(u,v)$ to $(gu,gv)$. Clearly, we have a homomorphism
$$
\Hie(\T)\to S(A\bigl| B)
$$
sending $\Aut(\T)$ to $S(A)\times S(B)$. Moreover, 
$\Aut(\T)$ is precisely the preimage of $S(A)\times S(B)$.
This implies the statement a).

On the other hand the image of $\Hie(\T)$ is closed 
in $S(A\bigl| B)$ and a closed subgroup of a Polish group
is Polish.

\sm 

{\bf \punct The ball-algebra $\cA$.}

\begin{lemma}
	Let $S_1$, $S_2$ be $(\T)$-subtrees of $\T$.
	
	\sm
	
	{\rm a)} If 
	$S_1\cap S_2$ is not empty, then
	 it is  a $(\T)$-subtree.
	
	\sm
	
		{\rm a)} If 
	$S_1\setminus S_2$ is not empty, then 
		$\vert(S_1)\setminus \vert(S_2)$ is a set of vertices of
		 a forest consisting of $(\T)$-subtrees.
\end{lemma}

This is obvious.

\sm 

Removing an edge of $\T$ we get two $(\T)$-subtrees. We call them 
 {\it branches}. We call a subset of the boundary $\partial\T$
 adjacent to a branch  a {\it ball}%
  \footnote{Let us use notation of Subsect..~\ref{ss:1.1}.
  	Consider a branch that do not cantain the initial point of
  	$\T$. Then the adjoint subset of the boundary is a ball in the sense
  	of the metric
 	$d_{\tochka}$.}.
 Consider the algebra%
 \footnote{We say that a family $\cA$ of subsets of a set $X$
 is an {\it algebra} if  $B\in \cA$ implies that
$X\setminus B\in\cA$ and $B$, $C\in\cA$ implies
$B\cap C$, $B\cup C\in\cA$.}
  $\cA(\T)$ of subsets in $\vert(\T)$
  generated by all branches. 
  By $\cA(\partial\T)$ denote the algebra of subsets
  in
  $\partial\T$ generated
 by all balls. 
 
 \begin{lemma}
 	{\rm a)}
 Any element of $\cA(\T)$ is a set of vertices of a 
 forest $S_1$, \dots, $S_k$ consisting of $(\T)$-subtrees.

 	{\rm b)}	The algebra $\cA(\T)$ is countable.
 	
 	\sm 
 	
 	{\rm c)}  The map sending $R\in \cA(\T)$ to its boundary
 	is an isomorphism of algebras $\cA(\T)$ and $\cA(\partial\T)$.
 	
 	\sm  
 	
 	{\rm d)} Sets $B\in \cA(\partial\T)$ are closed and open.

 	\sm
 	
 	{\rm e)} The group $\Hie(\T)$ acts on the set of nontrivial
 	elements of $\cA$ transitively. 
 \end{lemma}

{\sc Proof.}
Proofs of a)-d) are obvious. Let us prove e). Let $S_1$, \dots, $S_k\subset \T$
and $S_1'$, \dots, $S_{k'}\subset \T$ be two
 forests  of $(\T)$-subtrees. Denote by $T_1$, \dots, $T_m$ 
 and  $T_1'$, \dots, $T_{m'}'$ the complementary forests.
 We can subdivide any $(\T)$-subtree to several 
 $(\T)$-subtrees, therefore we can assume $k'=k$, $m'=m$.
 Now we take a spheromorphism sending $S_j\to S_j'$, $T_i\to T_i'$.
 \hfill 
 $\square$

\section{The Baire space and the Thompson group%
\label{s:thompson}}

\COUNTERS

{\bf \punct The correspondence between $\R\setminus \Q$
	and $\partial\T$.} Let $y_1$, $y_2\in\Q$.
Denote
$$
\ls y_1,y_2\rs:= (y_1,y_2)\cap \Q.
$$
Cut $\R\setminus \Q=\R\P^1\setminus \Q\P^1$ 
into 4 pieces $\ls-\infty,-1\rs$,
$\ls-1,0\rs$, $\ls0,1\rs$, $\ls1,\infty\rs$.
 We represent  points $x$ of these sets as
 continued fractions
\begin{align*}
\ls-\infty,-1\rs&: \quad x=-[s_0; s_1; s_2;\dots];
\\
\ls-1,0\rs&: \quad x=-[0; s_1; s_2;\dots];
\\
\ls0,1\rs&: \quad x=[0; s_1; s_2;\dots];
\\
\ls1,\infty\rs&: \quad x=[s_0; s_1; s_2;\dots].
\end{align*}

In all cases $s_j\in \N$. For definiteness,
consider $\ls0,1\rs$. We consider a tree 
whose vertices are enumerated by collections
$(s_1,\dots, s_k)$, edges have a form
\begin{equation}
\text{$(s_1,\dots, s_{k-1}$) --- $(s_1,\dots, s_{k-1},s_{k})$.}
\label{eq:---}
\end{equation}
We get a tree isomorphic to $\T$, and the boundary 
of this tree is identified with $\ls 0,1\rs$.

So we get 4 copies of the tree $\T$. Adding 3 edges
connecting their initial points
we unite them to one tree
 $\T$, the boundary of this
tree is in one-to-one correspondence  with $\R\setminus \Q$.
Thus we get the map
$$
\Xi:\R\setminus \Q \to \partial\T.
$$

Consider the algebra $\cR$ of subsets in $\R\setminus \Q$
generated by all intervals $\ls u,v\rs$ with rational
$u$, $v$ (we admit $u=-\infty$ and $v=\infty$).
On the other hand, we have the algebra
$$
\cA(\B):=\cA(\partial\T).
$$

\begin{proposition}
	The map $\Xi$ determines a bijection between algebras
	$\cR$ and $\cA(\B)$.
\end{proposition} 

{\sc Proof.} It is sufficient to show that any ball in
$\partial\T$ corresponds to an element of $\cR$ and
any interval $\ls u,v\rs$ corresponds to an element of
$\cA$.

\sm 

1) For definiteness let us remove an edge (\ref{eq:---})
in the $(\T)$-subtree corresponding to
$\ls 0,1\rs$. We get two branches of $\T$, one of them is completely
contained in the subtree. Its boundary consists of points
$$
[0;s_0,\dots,s_k, t_{k+1}, t_{k+2},\dots ], \quad \text{where $t_m$
	range in $\N$.}
$$
In other words, we get the interval
$$\LS[0;s_0,\dots,s_k],[0;s_0,\dots,s_{k}+1]\RS
\quad\text{or}\quad \LS[0;s_0,\dots,s_{k}+1],[0;s_0,\dots,s_k]\RS
$$
depending of the parity of $k$.

\sm  

2) Conversely, consider an interval $\ls u, v\rs$, where $0\le u<v\le 1$.
We have $\ls u, v\rs=\ls 0,w\rs\cap \ls v,	1\rs$.
For definiteness consider the subset of the Baire
space corresponding to $\ls 0,w\rs$.
Decompose $w$ into a continued fraction,
$$
w:=[0;s_1,\dots,s_k].
$$
Let $k=1$.
Then
$$
\ls 0,w\rs=\Ls 0, \tfrac{1}{s_1}\Rs
=
\ls 0, 1\rs\setminus \bigcup_{j=1}^{s_1-1} \Ls \tfrac 1{j+1}, \tfrac 1{j}  \Rs 
=
\ls 0, 1\rs\setminus \bigcup_{j=1}^{s_1-1} \Ls\, [0;j+1], [0;j]\,  \Rs.
$$
Let $k$ be even. We represent our interval as
\begin{multline*}
\Ls  0,[0,s_1,\dots,s_{k-1},s_k]\Rs=
\Ls  0,[0,s_1,\dots,s_{k-1}]\Rs
\\
\bigcup \bigcup_{j=1}^{s_{k}-1} \Ls\, [s_1,\dots,s_{k-1}+j], [s_1,\dots,s_{k-1}+j+1]
\,\Rs.
\end{multline*} 
For odd $k>1$ we write
\begin{multline*}
\Ls  0,[0,s_1,\dots,s_{k-1},s_k]\Rs=
\Ls  0,[0,s_1,\dots,s_{k-1}]\Rs
\\
\bigl\backslash \bigcup_{j=1}^{s_{k}-1} \Ls\, [s_1,\dots,s_{k-1}+j+1], [s_1,\dots,s_{k-1}+j]
\,\Rs.
\end{multline*} 
In all cases we have get an interval $\ls 0,w'\rs$,
where a continued fraction for $w'$ is shorter than
the a continued fraction for $w$, and a collection
of intervals $\ls p_j,q_j\rs$ corresponding to balls in $\partial \T$.
So we can apply the induction.
\hfill $\square$

\sm

{\bf \punct The action of the Thompson group.}

\begin{theorem}
	\label{th:thompson}
	Let $\kappa$ be a transformation  of $\R\setminus \Q$ 
	lying in the Thompson group $\Th$. Then $\Xi\circ \kappa \circ \Xi^{-1}$
	is contained in $\Hie(\T)$.
\end{theorem}

Recall that the group $\PGL_2(\Z)$ acts on
$\R\setminus\Q$ by linear fractional transformations.
First, we prove the following preliminary statement.

\begin{proposition}
	\label{pr:pgl}
For any	$h\in \PGL_2(\Z)$ we have
 $\Xi\circ h \circ \Xi^{-1}\in \Hie(\T)$.
\end{proposition}

{\sc Proof.} It is sufficient to prove this statement for
generators 
$$
\begin{pmatrix}
0&1\\1&0
\end{pmatrix}: \quad x\mapsto \frac 1x;
\qquad 
\begin{pmatrix}
1&1\\0&1
\end{pmatrix}: \quad x\mapsto  x+1
$$	
of the group $\PGL_2(\Z)$.

\sm  

a) Let $x\in \ls 0,1\rs$,
$$
x=[0;s_1,s_2,s_3,\dots].
$$
Then
$$
x^{-1}=[s_1;s_2,s_3,\dots]
$$
and we  permute  branches $\Xi\ls 0,1\rs$
and
$\Xi\ls 1,\infty\rs$ preserving their structures.
The same holds for  $\ls-1,0\rs$ and $\ls-\infty,-1\rs$.

\sm 

b) Examine the transformation $x\mapsto x+1$.

\sm 

b.1) Let $x\in \ls -\infty,-2\rs$, 
$$
x=-[s_0;s_1,s_2,\dots],\quad x+1=-[s_0-1;s_1,s_2,\dots].
$$
Transformation of subtrees corresponding 
to the shift
$\ls-\infty,-2\rs\to \ls-\infty,-1\rs$
is shown on 
 Fig. \ref{fig:shift1}.
\begin{figure}
$$\epsfbox{inftree.5}$$	
\caption{Ref. to the proof of Proposition \ref{pr:pgl}. The shift $\ls(-\infty,-2\rs\to \ls(-\infty,-1\rs$.}
\label{fig:shift1}
\end{figure}
%

\sm 

b.2) The map $\ls -2,-1\rs\to\ls -1,0\rs$.
We send $x=[-1;s_1,s_2,\dots]$ to $x=[0;s_1,s_2,\dots]$.
This is an isomorphic map of two branches of $\T$.

\sm

b.3) Let $x\in \ls -1,0\rs$. We represent it in two forms
\begin{align*}
x=-[0;s_1,s_2,\dots]&=-\frac 1{s_1+\eta},\qquad \text{where $0<\eta<1$},
\\
 &=-\cfrac{1}{s_1+\cfrac{1}{s_2+\xi}},
\qquad \text{where $0<\xi<1$.}
\end{align*}
Let $s_1>1$. Then
$$
-x+1=-\frac 1{s_1+\eta}+1=\frac{s_1-1+\eta}{1+\eta}=
\cfrac{1}{1+\cfrac{1}{s_1-1+\eta}}.
$$

For $s_1=1$ we have
$$
-x+1=-\cfrac{1}{1+\cfrac{1}{s_2+\xi}}+1
=\cfrac{1}{s_2+1+\xi}.
$$
The correspondence of branches is shown on Fig. \ref{fig:shift2}.
\begin{figure}
	$$ \epsfbox{inftree.6} $$
	\caption{Ref. to the proof of Proposition \ref{pr:pgl}. The shift $\ls-1,0\rs \to \ls 0,1\rs$.}
	\label{fig:shift2}
\end{figure}
%


\sm 

b.4) The examination of the shift $x\mapsto x+1$
on $\ls 0,1\rs $, $\ls 1,\infty\rs$ is similar to the case b.1).
\hfill $\square$

\sm 

{\sc Proof of Theorem \ref{th:thompson}.} 
Consider a piece-wise $\PSL_2(\Z)$-transformation $h$ of $\R\P^1$.
The projective line is a union of of rational
segments $[r_m,r_{m+1}]$, on which the transformation
corresponds to some elements $\gamma_m\in \PSL_2(\Z)$.
Such segments $[r_m,r_{m+1}]$ correspond to elements $B_m$ of the algebra $\cA$,
for each element we have a finite collection $S_{m1}$, $S_{m2}$,
\dots  of
$(\T)$-subtrees. On the other hand $\gamma_m$ determines
a spheromorphism and therefore a finite collection of
$R_{m1}$, $R_{m2}$, $\dots$ of $(\T)$-subtrees.
Therefore
$$
U_{mij}:= S_{mi}\cap R_{mj}
$$
is a splitting of $\T$ into $(\T)$-subtrees and $\Xi \circ h\circ \Xi^{-1}$ embeds
each subtree to $\T$.
\hfill $\square$

\section{The train of the group $\Hie(\T)$%
\label{s:train}}

\COUNTERS

{\bf \punct Double cosets and combinatorial data.%
	\label{ss:bi-tree}}
For a {\it non-empty} finite subtree $J\subset \T$
denote by $\cK(J)\subset \Aut(\T)$ the pointwise stabilizer of $J$.
We wish to describe  spaces of double cosets
$$
\cK(J)\setminus \Hie(\T)/\cK(I)$$
for finite nonempty subtrees $I$, $J$.

Consider a graph $\Gamma$ whose edges are colored black, blue, and red.
Denote by $\Gamma^\black$ the subgraph consisting of black edges,
by $\Gamma^\black_\blue$ of black and blue edges,
by $\Gamma^\black_\red$ of black and red edges.

Let $\Gamma$ be such a graph equipped with 
embeddings $\imath:I\to \Gamma$, $\jmath:J\to\Gamma$.
We say that $(\Gamma,\imath,\jmath)$ is a {\it $(I,J)$-bi-tree} 
if the following conditions hold:

\sm

$\bullet$
the subgraphs $\Gamma^\black_\blue$ and
$\Gamma^\black_\red$ are trees;

\sm

$\bullet$ $\imath(I)\subset \Gamma^\black_\blue$,
$\jmath(J)\subset \Gamma^\black_\red$;

\sm

$\bullet$
any vertex of $\Gamma$ of valence 1 is an end of a black
edge contained in $\imath(I)$ or $\imath(J)$.

\sm 

\begin{figure}
$$\epsfbox{inftree.7}$$

Notation:
 $\epsfbox{segment.1}$
  is black,
$\epsfbox{segment.3}$
 is blue, 
$\epsfbox{segment.2}$ 
is red,
\\
edges
$\epsfbox{segment.4}$ 
are in $\imath(I)$,
$\epsfbox{segment.5}$
 are in $\jmath(J)$.

\caption{Ref. to Subsect. \ref{ss:bi-tree}.
An $(I,J)$-bi-tree. Removing blue and red edges
we get a black forest consisting of 4 trees.}
\label{fig:bi-tree}
\end{figure}

Two $(I,J)$-bi-trees $(\Gamma,\imath,\jmath)$
and $(\Gamma',\imath',\jmath')$ are {\it equivalent}
if there is a color preserving isomorphism 
$\lambda: \Gamma\to\Gamma'$ such that $\imath'=\lambda \circ \imath$,
$\jmath'=\lambda \circ \jmath$.

\sm 


\begin{proposition}
	There is a canonical one-to-one correspondence
	between the space $\cK(I)\setminus \Hie(\T)/\cK(J)$
	of double cosets and the set
	of all $(I,J)$-bi-trees defined upto the equivalence.
\end{proposition} 

This is version of the correspondence defined in Subsect. \ref{ss:bi-tree1}.

\sm 

{\sc Construction of an $(I,J)$-bi-tree by a spheromorphism
	$g$.}
Fix $g\in \Hie(\T)$. Consider the perfect
 $(\T)$-covering forest $\{S_\alpha \}$ compatible with $g$,
 see Subsect. \ref{ss:canonical-forest}.
 Paint blue all edges in $\T\setminus \bigcup S_j$.
 The {\it $I$-skeleton $\Skel_I\{S_\alpha \}$} of $\{S_\alpha \}$ is the minimal
 subgraph in $\T$ containing all blue edges and $I$.
 Next, we paint red all edges of $\T\setminus \bigcup gS_j$,
 and consider $J$-skeleton of $\{gS_\alpha \}$,
 i.~e., the minimal subgraph in $\T$ containing red edges
 and $J$.
 
 In each $S_k$ we have two subtrees, $S_k\cap\Skel_I(\{S_\alpha \})$
 and $g^{-1} (gS_k \cap \Skel_J(\{gS_\alpha \}))$.
Consider the minimal  subtree $\Xi_k$ containing these subtrees
and paint it  black. Any blue edge in $\Skel_I(\{S_\alpha \})$
has two ends in some $S_k$, $S_m$. These ends are
contained in $\Xi_k$, $\Xi_m$.
So we add blue edges to the forest $\{\Xi_\alpha \}$,
in the same way we add red edges.

Since our $(\T)$-covering forest is perfect,
we do not get double red-blue edges.
 
\sm 

{\sc Remark.} For a given spheromorphism $g$ 
we have a canonical embedding $i_+$ of $\Gamma^\black_\blue$
to $\T$ and a canonical embedding $i_-$ of $\Gamma^\black_\red$
to $\T$. They are related by 
$$
i_-\Bigr|_{\Gamma^\black}=g\circ i_+\Bigr|_{\Gamma^\black}
$$
We denote these subtrees in $\T$ by $^*\Gamma^\black_\blue$ and $^*\Gamma^\black_\red$.
\hfill $\square$

\sm 

{\sc The inverse construction.}
Let $\Gamma$ be an $(I,J)$-bi-tree. 
Take two copies $\T_1$, $\T_2$
of the tree $\T$ with $I$ drawn on $\T_1$ and $J$
drawn on $\T_2$. Choose isomorphisms $\theta_{1,2}:\T\to \T_{1,2}$.
Consider  embeddings 
$$p: \Gamma^\black_\blue\to \T_1,\qquad q: \Gamma^\black_\red\to \T_2$$
such that
$$
\theta^{-1}_1\circ p\circ \imath,\quad 
\theta^{-1}_2\circ q\circ \jmath
$$
are identical maps $I\to I$ and $J\to J$ respectively.

Then $p$  determines a $(\T)$-covering forest, say $\{P_\alpha \}$,
of $\T_1$ (we remove images of blue edges from
$\T_1$).
In the same way $q$ determines a   $(\T)$-covering forest,
say $\{Q_\alpha \}$, of $\T_2$. 

The set of components of $\{P_\alpha \}$ (resp., $\{Q_\alpha \}$)  is in
a one-to-one correspondence with the set
of components $\{\Xi_\alpha\}$
of the black forest $\Gamma^\black$.
 For each 
$P_k$ we restrict $p$ to $\Gamma^\black_\blue\cap \Xi_k$
and get an embedding of this tree to $P_k$. Extend it
to an embedding $\wt p_k:\Xi_j\to P_k$.

In the same way we get embeddings $\wt q_j:\Xi_j\to Q_j$.
Next, we choose isomorphisms $\wh r_j:P_j\to Q_j$
such that 
$$
\wt q_j\Bigr|_{\Xi_j}=\wh r_j\circ \wt p_j\Bigr|_{\Xi_j}
.$$
This determines a spheromorphism
$r:=\{\wh r_j\}:\T_1\to \T_2$
and  a spheromorphism
$g:\T\to\T$ 
$$
g=\theta_2^{-1}\circ r\circ \theta_1.
$$

Multiplying $\theta_1\mapsto \theta_1 h_1$,
$\theta_2\mapsto\theta_2 h_2$, where
$h_1\in \cK(I)$, $h_2\in \cK(J)$,
we get all elements of the double coset.

On the other hand maps $\wh r_k$ are not canonical
and they can be replaced by maps $\wh r_k\circ  \zeta_k$,
where $\zeta_k$ are maps $P_k\to P_k$ fixing $\Xi_k$.
This determines a spheromorphism $\{\zeta_k \}$,
which is contained in $\Aut(\T)$. 
Now we can replace $\theta_1$
by $\{\zeta_k \}\circ\theta_1=\theta_1\circ (\theta_1^{-1}\circ \{\zeta_k \}\circ\theta_1)$
and we get an element of the same double coset. 


{\sc Weak bi-trees of spheromorphisms.}
In proof of Lemma \ref{l:gg} we need a variation of the construction.
Consider the following data: a spheromorphism $g$, a compatible $(\T)$-covering
forest $\{S_\alpha \}$ (generally, non-perfect)
and a collection $\{v_j\}$ of marked vertices in $\T$. Then we  apply
the procedure of drawing of  a bi-tree and get the graph
$\varGamma$ whose edges are colored black, blue, red, and 
double blue-red edges are allowed. We define
the subgraph $\varGamma^\black_\blue$ as the graph obtained
by  removing red edges (double blue-red edge become blue).
In a similar way we define the subgraph $\varGamma^\black_\red$.
Again, $\varGamma^\black_\blue$, $\varGamma^\black_\red$
are trees, whose terminal black edges finish at marked points. 

\begin{figure}
	$$\epsfbox{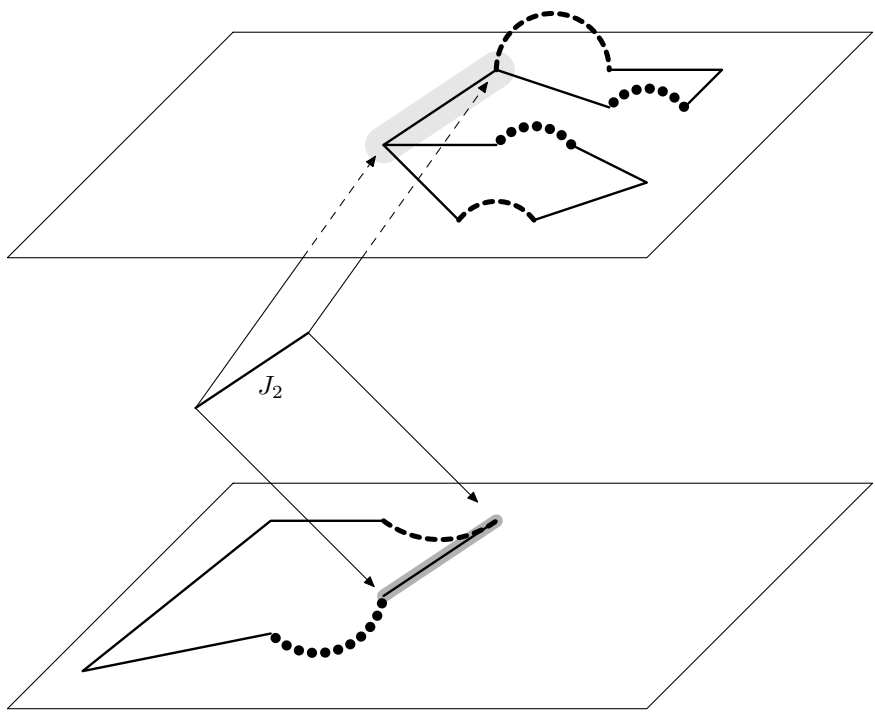}$$
	\caption{Refs. to 
		Subsect. \ref{ss:category} and Lemma \ref{l:independent}. Gluing of bi-trees.}
	\label{fig:product}
\end{figure}

{\bf\punct The category of bi-trees and the category
	of double cosets.%
	\label{ss:category}}
Denote by $\cM(I,J)$ the set of $(I,J)$-bi-trees.
We wish to define a product 
$$
\cM(J_1,J_2)\times \cM(J_2,J_3)
\to \cM(J_1,J_3).
$$
Let $\Gamma\in\cM(J_2,J_3)$, $\Delta\in\cM(J_1,J_2)$.
We glue $\Gamma$ and $\Delta$ identifying
images of embeddings $J_2\to \Gamma$,  $J_2\to \Delta$, see Fig.\ref{fig:product}.
After this we can get double colors on some edges of
$J_2$. We replace%
\footnote{We indicate both color of an edge and the origin of an edge
	($\Delta$ or $\Gamma$).}
$$
\text{($\Delta$-blue, $\Gamma$-black)}\longrightarrow \text{(blue)},\quad
\text{( $\Delta$-black, $\Gamma$-red)}\longrightarrow \text{(red)}$$
and remove
($\Delta$-blue, $\Gamma$-red)-edges%
\footnote{We also have ($\Delta$-black, $\Gamma$-black)$\longrightarrow$
	(black), other combinations of colors are impossible.}.
Finally, remove all vertices of valence 1,
that are not contained in the images of $J_1$ and $J_3$
(such vertices are automatically contained in $J_2$
and adjacent edges are black). Repeat the step again, etc%
\footnote{A black edge survives if and only if
	it can be included to a way $v_1,\dots,v_m$ following black edges,
	whose terminal vertices $v_1$, $v_m$ are contained in $J_1$, $J_3$ or are ends of
	blue or red edges.}.

\begin{lemma}
	In this way, we get a $(J_1,J_3)$-bi-tree.
\end{lemma}

Denote this $(J_1,J_3)$-bi-tree by
$\Delta \diamond \Gamma$.

\sm

{\sc Proof.} Let us examine the graph $\Theta$
obtained by glueing of $\Delta$ and $\Gamma$.
Consider its subgraph consisting of edges, on which black or
red are present%
\footnote{Recall that some edges have two colors.},
i.e., $\Delta^\black_\red\cup \Gamma^\black_\red$.
The subtree $\Gamma^\black_\red$ contains
$J_2$. The graph $\Delta^\black_\red\setminus\{\edge(J_2)\}$ is a forest
and each component has one or two vertices in $\vert(J_2)$.
Consider an edge $[v,w]$ of $J_2$. There are two cases:

\sm 

1) The edge $[v,w]$ is black in $\Delta$. This edge is a unique way
in $\Delta^\black_\red$
connecting $v$ and $w$. Since $\Gamma^\black_\red\supset J_2$ can be contracted
to $J_2$, then   $[v,w]$ is a unique way in $ \Delta^\black_\red\cup\Gamma^\black_\red$
connecting $v$ and $w$.

\sm 

2) The edge $[v,w]$ is blue in $\Delta$. According our rules
it is blue in $\Delta\diamond \Gamma$ and absent in
$[\Delta\diamond \Gamma]^\black_\red$. However, the vertices $v$ and $w$
are connected by a unique
way in $\Delta^\black_\red$, and  this way is contained in
$[\Delta\diamond \Gamma]^\black_\red$.

\sm 

The same argument holds for $[ \Delta\diamond\Gamma]^\black_\blue$.

Next, all edges of $J_2\cap J_3$ and of $J_1\cap J_2$
are black. Therefore edges of $J_1$ and $J_3$ can not disappear
after removing  blue-red edges.
\hfill $\square$

\begin{lemma}
	Let $\Lambda$ be a $(J_1,J_2)$-bi-tree, $\Delta$ a $(J_2,J_3)$-bi-tree,
	$\Gamma$ a $(J_3,J_4)$-bi-tree. Then
	$(\Lambda\diamond\Delta)\diamond \Gamma =\Lambda\diamond(\Delta\diamond \Gamma) $.
\end{lemma}

{\sc Proof.} We glue  $\Lambda$, $\Delta$,  $\Gamma$
identifying two copies of $J_3$ in  $\Delta$ and $\Gamma$,
and $J_2$ in $\Lambda$  and $\Delta$. Clearly, order of gluings has no matter.
Different orders of recolorings can appear, when $J_2$ and $J_3$  have a common edge.
This edge must be black in $\Delta$.
In $\Gamma$ it can be red or black, in $\Lambda$ black or blue.
In all admissible four cases result does not depend on order of recolorings.
\hfill $\square$

\sm 

{\sc Remark.} If the tree $J_2$ is empty
then this product is not well-defined (since we do not get a tree).
\hfill $\square$

\sm 

Thus we get a category 
whose objects are (non-empty) finite
trees, and  morphisms $I\to J$ are $(I,J)$-bi-trees. 

Since $\cM(J,I)$ is in one-to-one correspondence with
double cosets, we get
the product of double cosets 
$$
\cK(J_1)\setminus \Hie(\T)/\cK(J_2)
\,\times \, \cK(J_2)\setminus \Hie(\T)/\cK(J_3)
\,
\to\,
\cK(J_1)\setminus \Hie(\T)/\cK(J_3).
$$
Denote this category $\hier$. Objects are 
{\it nonempty} finite subtrees $J$ in $\T$,
morphisms are 
$$
\Mor_\hier(J_2,J_1):=\cK(J_1)\setminus \Hie(\T)/\cK(J_2)\simeq \cM(J_1,J_2).
$$
We denote the multiplication of morphisms in $\hier$ by $\odot$.


{\bf \punct Lemma about independence.}
\begin{lemma}
	\label{l:independent}
	Let $J_1$, $J_2$, $J_3\subset\T$
	be nonempty finite subtrees.
	Let $g_1$, $g_2$ be spheromorphisms, 
denote by	$\Delta$  the $(J_2,J_1)$-bi-tree of the spheromorphism
	 $g_1$,
by	$\Gamma$  the $(J_3,J_2)$-bi-tree of  $g_2$.
	Let $\Xi$ be the   $(J_3,J_1)$-bi-tree of the product
	 $g_1 g_2$.
Assume that	
	$$
	^*\!\Gamma^\black_\red\cap\, ^*\!\Delta^\black_\blue=J_2.
	$$
	Then $\Xi=\Delta\diamond\Gamma$.
\end{lemma}


{\bf \punct Bi-trees of products. Proof of Lemma \ref{l:independent}.%
	\label{ss:of-products}}
Let $g_1$, $g_2\in\Hie(\T)$ be arbitrary.
We intend to describe the $(J_3,J_1)$-bi-tree $\Theta$ 
of $g_1g_2$ if we know $(J_1,J_2)$-be-tree $\Delta$ of $g_1$, 
$(J_2,J_3)$-bi-tree $\Gamma$ of $g_2$, and the map $g_2^{-1}$ 
on $\vert(^*\!\Delta^\black_\blue)$.

Consider the union of the subtrees 
$$\Theta^{\circ\circ}:=^*\!\!\Delta^\black_\blue\cup\, ^*\!\Gamma^\black_\red\subset \T,$$
color edges of $\T$ outside this union as {\it grey}.
The intersection of these subtrees contains $J_2$ and hence it is not empty.
Therefore $\Theta^{\circ\circ}$ is a subtree. So we colored $\T$ in 4 colors,
grey, black, blue, red (some  edges are colored in two colors).
We add to $\Theta^{\circ\circ}$ red edges of $\Delta$ and blue edges 
of $\Gamma$ and get a new graph $\Theta^{\circ}\supset \T$, some its edges
can be double.  
We consider this picture as a pair (graph $\Theta^\circ$, subgraph $\T$).
So we can distinguish red  edges originated from
$\Gamma$ (they are contained in $\T$) and from $\Delta$
(they are not contained in $\T$). 

We claim, that the bi-tree $\Theta$ of $g_1g_2$ is obtained 
from $\Theta^\circ$ by the following way:

\sm 

\begin{figure}
		$$\epsfbox{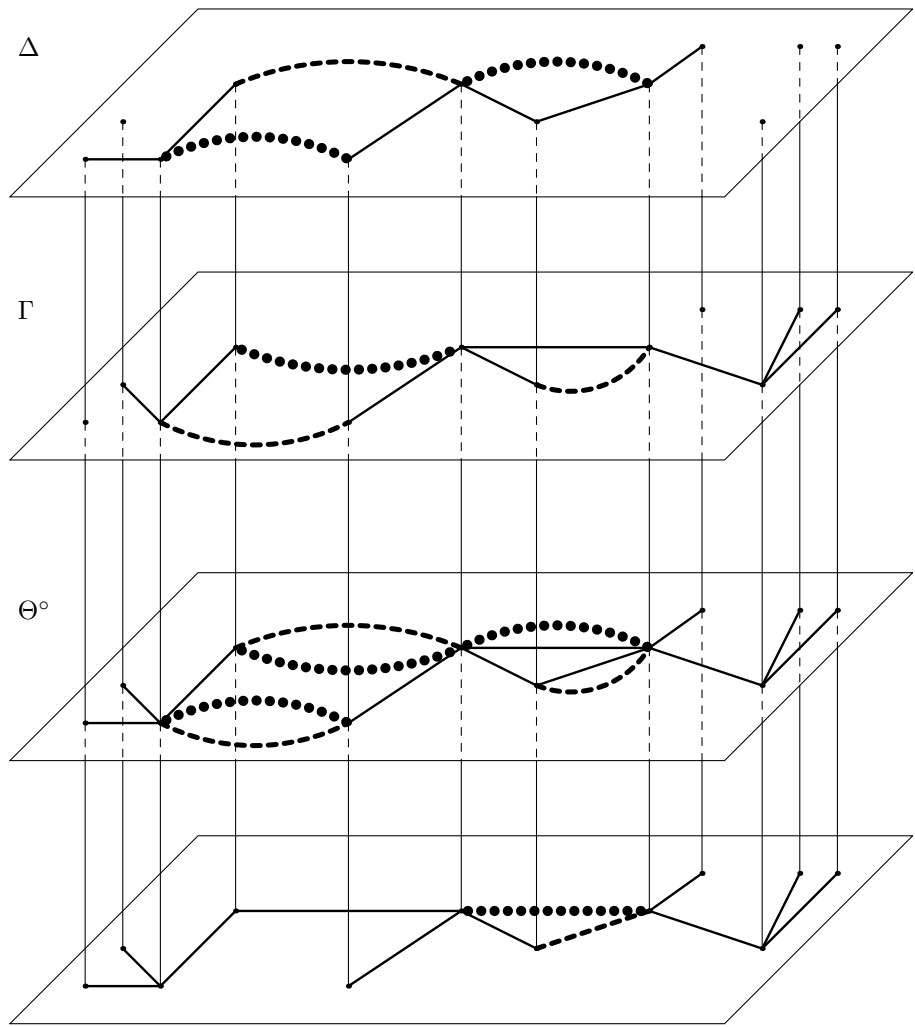}$$
	\caption{Ref. to Subsect. \ref{ss:of-products}.
		The recoloring. We draw pieces of $\Delta$, $\Gamma$, the corresponding
		piece of $\Theta^\circ$ (grey edges are omitted), and the result of application 
		of steps $A^*$--$B^*$ of the recoloring.}
	\label{fig:recoloring}
\end{figure}

$A^*)$ we transform double edges to simple edges according their colorings,
\begin{gather*}
\text{($\Delta$-blue,$\Gamma$-black)}\longrightarrow \text{($\Theta$-blue)},
\qquad
\text{( $\Delta$-black, $\Gamma$-red)}\longrightarrow \text{($\Theta$-red)};
\\
\text{( $\Delta$-red, $\Gamma$-blue)}\longrightarrow \text{($\Theta$-black)};
\end{gather*}

$B^*)$  remove double edges of the type
($\Delta$-blue, $\Gamma$-red);

\sm

$C^*)$ remove grey edges;

\sm 

$D^*)$  in the rest we successively remove all terminal black
edges that are not contained in the images of $J_1$ and $J_3$.

\sm 

See Fig. \ref{fig:recoloring}.

Keeping in mind the  proof of Lemma \ref{l:gg} below, we
denote by $\Theta^{\max}$ the result of application
of operations $A^*$--$C^*$ to $\Theta^\circ$.


\begin{lemma}
	\label{l:pr}	
	The graph  $\Theta$ is a $(J_1,J_3)$-bi-tree.	
\end{lemma}

{\sc Proof.} Examine the transformation of $\Delta\subset \Theta^\circ$
under changing of colors. The subtree $\Delta^\black_\red\subset \Delta$
remains to be colored black and red, but some  black edges can became
red and some red edges can became black. On the other hand blue edges
of $\Delta$ can disappear but they can not be recolored black or red.
So no black or red edges can be added to $\Delta^\black_\red$.
Thus edges that are contained simultaneously in $\Theta^\black_\red$ and $\Delta$
form a tree.

On the other hand, $\Gamma^\black_\red$ is a subtree in $\T$.
Remove edges that are contained in $\Delta$, 
 $$
\edge(\Gamma^\black_\red)\setminus
\edge(\Delta)=\edge(\Gamma^\black_\red)\setminus \edge(\Delta_\black^\blue).
$$
We get a difference of two subtrees in $\T$, it is a forest.
 Each component of this forest has a unique vertex common with
$\Delta^\black_\red$. So  $\Theta^\black_\red$ is a tree.

The same argument shows that $\Delta^\black_\blue$
is a tree.

Next, the image of $J_1$ in $\Delta$  consist of black and 
red edges. As we have seen these edges can be recolored 
but can not disappear.

Thus, after application of the transformations $A^*$--$C^*$
to $\Theta^\circ$ we get a graph satisfying all properties of
$(J_1,J_3)$ except the absence of terminal edges.
Such edges disappear after the transformation $D^*$.
\hfill $\square$


\begin{lemma}
	\label{l:gg}
	The graph $\Theta$ is the $(J_1,J_3)$-bi-tree of 
	$g_1g_2$.
\end{lemma}

{\sc Proof.} Let $\T_1$, $\T_2$ be copies of $\T$. Let us think
that $g_2$ sends $\T\to\T_1$ and $g_1:\T_1\to\T_2$.
Denote by $\Theta^\lozenge$ the $(J_1,J_3)$-bi-tree of $g_1g_2$.

\sm

{\it An upper estimate of $\Theta^\lozenge$.}
On  $\T$ we mark some edges and vertices according
the following rules.
Consider the perfect $(\T)$-covering forest $\{S_\alpha\}$ for $g_2$ and paint
 a color {\it $\Gamma$-blue} all edges in $\T\setminus \cup S_\alpha$.
Also paint  $\Gamma$-blue their ends. 
Paint edges of $\T_1\setminus \cup g_2S_\alpha$ to a color {\it $\Gamma$-red},
also paint $\Gamma$-red their ends on $\T_1$ and $g_2$-preimages of their ends
in $\T$ (we admit several colors at one vertex).

Next, take the perfect $(\T)$-covering forest for $g_1$ drawn on $\T_1$.
Repeat the same procedure with colors {\it $\Delta$-blue} and {\it $\Delta$-red}.
On the initial copy $\T$
we paint  $\Delta$-blue $g_2$-preimages of $\Delta$-blue
edges (if a preimage of an edge
is an edge)
and preimages of $\Delta$-blue vertices.
We also paint $\Delta$-red $g_1g_2$-preimages of $\Delta$-red vertices.

Finally, we mark points of the sets 
$$\text{$\vert(J_3)$, $g_2^{-1}\vert(J_2)$,
$g_2^{-1}g_1^{-1}\vert(J_1)\subset \T$,}$$
 we call such vertices $J_3$-vertices,
$J_2$-vertices, $J_1$-vertices.

Consider the minimal subtree $\Sigma$ in $\T$ containing all marked data.
Clearly%
\footnote{Generally, the inclusion is strict, generally 
	$\vert(\Theta^\lozenge)$ does not contain $g_2^{-1}\vert(J_2)$;
	also  $g_1$ can restore an edge cut by $g_2$ but on our picture this 
	event leaves marked vertices.},
$$\Sigma\supset [\Theta^\lozenge]^\black_\blue
\qquad\text{(as  non-colored graphs)}.$$

{\it A description of $\vert(\Sigma)$ in the terms of $\Gamma$ and $\Delta$}.
Removing $\Delta$-blue and $\Gamma$-blue  
 edges from $\T$ we get a $(\T)$-covering
forest
for  $g_1g_2$, say $\{Q_\nu \}$ (it can be non-perfect). Denote by $\Sigma^\dagger$
the forest obtained from $\Sigma$ by removing $\Delta$-blue and $\Gamma$-blue
edges.

Each  $\Sigma^\dagger_j$ is a minimal subtree in the corresponding
$(\T)$-subtree $Q_j$ containing all marked vertices,
i.e., vertices of the types
\begin{equation}
\text{ $\Delta$-blue, $\Gamma$-blue,
	$\Delta$-red, $\Gamma$-red, $J_{1}$, $J_2$, $J_3$}.
\label{eq:set1}
\end{equation}

Examine the corresponding black subtrees in $\Delta$ and $\Gamma$.

For $g_2$, $(\T)$-covering forest $\{Q_\nu\})$,  marked
$J_{2}$-vertices and $J_3$-vertices we construct
the corresponding weak bi-tree%
\footnote{If we replace double $\Gamma$-blue--$\Gamma$-red
	edges of $\varGamma$ by single black edges we get the same
	graph $\Gamma$.}
$\varGamma$ as at the end of Subsect. \ref{ss:bi-tree}.
Consider the corresponding black forest
$\{\varXi_\nu\}$.
Its element $\varXi_j$
is the minimal subtree
in $Q_j$ containing all
vertices of the following types
\begin{equation}
\text{ $\Gamma$-blue,  $\Gamma$-red,
	$J_2$, $J_3$}.
\label{eq:set2}
\end{equation}
For $g_1$, the $(\T)$-covering forest $g_2\{Q_\nu \}$, and 
marked set $J_2\cup g_1^{-1}J_1$ consider the corresponding
weak bi-tree $\varDelta$. Consider the corresponding black forest
and its $g_2$-preimage $\{Z_\nu\}$. 
The tree $Z_j$  is the minimal subtree
in $Q_j$ containing  all vertices
of the following types:
\begin{equation}
\text{$\Delta$-blue, $\Delta$-red,
	$J_1$, $J_2$}.
\label{eq:set3}
\end{equation}
Since (\ref{eq:set1}) is a union of (\ref{eq:set2}) and (\ref{eq:set3}),
we 
come to the following alternative:
$$
\Sigma^\dagger_j=\varXi_j\cup Z_j\qquad \text{or}\qquad \varXi_j\cup Z_j=\varnothing.
$$


\begin{lemma}
	\label{l:ischo}
	We have $\Sigma^\dagger_j=\varXi_j\cup Z_j\qquad$.
\end{lemma}

This implies that
$$
\vert(\Sigma^\dagger)\subset \vert(\Gamma)\cup g_2^{-1}\vert(\Delta):=
\vert(^*\!\Gamma^\black_\blue)\cup g_2^{-1}\vert(^*\!\Delta^\black_\blue).
$$

{\sc Proof of Lemma \ref{l:ischo}.}
Assume the contrary. Let $x\in \varXi_j$ and $z\in Z_j$ be the nearest vertices of $\varXi_j$ and $Z_j$.
Let $[x, v_1,\dots, v_m,z]$ be the way connecting $x$ and $z$.
Cutting the first and the last edges of this way we get a $(\T)$-covering forest,
consisting of three or two pieces, $A$ containing $x$, $C$ containing $z$, 
and the rest $B$, which can be empty. 

Then there are no  vertices of types (\ref{eq:set2})
in $B\cup C$. Otherwise 
there is a way on $\T$ connecting $\varXi_j$ with such a vertex, 
the first edge of the way,
namely $[x,v_1]$, must be black and therefore must be contained in $\varXi_j$.

Also there are no vertices of types (\ref{eq:set3})
in $A\cup B$. Indeed, consider a way on $\T_1$ connecting $g_2 Z_j$ with such a vertex.
Its preimage on $\T$ is a collection of black ways whose ends are $\Gamma$-red.
However, such a 'path' can not leave $C$, indeed there no $\Gamma$-red vertices
in $C$, so a jump is impossible, on the other hand a continuous way can not avoid
the edge $[v_m,z]$ but it is not black.

Thus $J_2$-vertices can not be contained in $A$, $B$, $C$.
We get a contradiction.
\hfill $\square$   









\sm 

{\sc End of proof of Lemma \ref{l:gg}.} 
{\it Rules of the recoloring.}
Next, we must examine the actual presence of edges
in the bi-tree  $\Theta^\lozenge$ and their colors.

\sm 

1$^*$. Consider a double edge $[v,w]$ of the type
($\Delta$-blue, $\Gamma$-black).
This means that we have two vertices $v$, $w$ such that
$[v,w]$ is an edge in $\T$, $[g_2v,g_2w]$ also is an edge
and $g_1g_2v$, $g_1g_2w$ are not connected by  an edge
in $\T$. So our edge of the bi-tree $\Theta^\lozenge$ is blue.

\sm 

2$^*$. The  similar argument holds for the combination 
($\Delta$-black, $\Gamma$-red).

\sm

3$^*$. Consider an edge of the type ($\Delta$-blue, $\Gamma$-red).
We have a pair of vertices $v$, $w\in \T$, which are not connected by
an edge, the edge $[g_2v,g_2w]$, and $g_1g_2v$, $g_1g_2w$,
which
are not connected by an edge in $\T$.
So we have no corresponding edge in $\Theta^\lozenge$.

\sm  

4$^*$. Consider a double edge of $\Theta^\circ$ of the type
($\Delta$-red, $\Gamma$-blue)
\footnote{Notice that both copies of the edge are
	not contained in $\T$.}.
This means that we have two vertices
$v$, $w$ of $\T$ such that $[v,w]$ is an edge of $\T$,
$g_2v$ and $g_2w$ are not connected by an an edge, and $[g_1g_2v,g_1g_2w]$
is again an edge. Therefore $[v,w]$ is not blue and 
$[g_1g_2v,g_1g_2w]$ is not red. So if this edge is present in
$\Theta^\lozenge$, then it is black.
Paint it {\it yellow}. 

\sm 

{\it The graphs $\Theta^{\max}$ and $\Theta^\lozenge$.}
Consider the forest $\{\Sigma^\dagger_\nu\}$.
Notice that for each vertex of the types $\Delta$-blue, $\Gamma$-blue,
$\Delta$-red, $\Gamma$-red
in $\Sigma^\dagger_j$  there is a corresponding vertex
of the same type in 
another tree $\Sigma^\dagger_j$ (since each colored
vertex appeared as an end of a colored edge). We draw the corresponding edges,
recolor  the graph as above, paint yellow edges to black
and get a new graph. It is clear that it is 
the graph $\Theta^{\max}$ defined above in this subsection.
By construction, $\Theta^{\max}\supset \Theta^\lozenge$.

We know the perfect $(\T)$-covering forest for
$g_1g_2$. Namely, if $\Sigma_k^\dagger$ and $\Sigma_m^\dagger$ are
connected by a yellow edge, then we connect
$(\T)$-subtrees $S_k$ and $S_m$ by an edge and unite them
to one $(\T)$-subtree. So we can describe $\Theta^\lozenge$.

Define the following  set of distinguished vertices of $\Theta^{\max}$:
\begin{equation}
\text{ends of blue edges, end of red edges, vertices originated from $J_1$ or $J_3$.}
\label{eq:set-new}
\end{equation}
Now we can formulate the following criterion:

\sm 

--- a black edge $[a,b]\in \Theta^{\max}$ is contained $\Theta^\lozenge$
if and only if it can be included to a way $v_1$, \dots, $v_N$,
consisting of black (or yellow) edges and
the ends $v_1$, $v_N$ of the way are contained in the set (\ref{eq:set-new}).

\sm

The subgraph $\Theta\subset \Theta^{\max}$
is a $(J_1,J_3)$-bi-tree, so its black edges satisfy this criterion,
therefore $\Theta\subset \Theta^\lozenge$.
On the other hand, $\Theta^{\max}\setminus \Theta$ is a forest.
Each its component has one vertex in $\Theta$, the remaining vertices
are not distinguished and therefore  edges of the component are not contained 
in $\Theta^\lozenge$
\hfill $\square$ 

\sm

{\sc Proof of Lemma \ref{l:independent}.}
We evaluate a bi-tree of the product according the prescription.
\hfill $\square$


\sm 

{\bf\punct Representations of the category of double cosets.}
  Let $\rho$ be a unitary representation of the group
  $\Hie(\T)$ in a Hilbert space $H$. For a finite subtree
  $J\subset \T$ denote by $H(J)$ the subspace of 
  $\cK(J)$-fixed vectors.
  If $J_1\subset J_2$, then $\cK(J_1)\supset \cK(J_2)$ and $H(J_1)\subset H(J_2)$.
   By $P(J)$ we denote the operator
  of orthogonal projection to $H(J)$.
  
  \begin{lemma}
  The subspace $\cup_J H(J)$
  is dense in $H$.	
  \end{lemma} 
  
  This is a special case of the following statement, see  см. \cite{Ner-book},
  Proposition VIII.1.2.
  
  \begin{proposition}
  	\label{pr:ism}
  	Let $G$ be a topological group, $Q_1\supset Q_2\supset \dots$
  be a family of subgroups such that each neighborhood
  of the unit in $G$ contains a subgroup $Q_j$. Consider a unitary 
  representation of $G$ in a Hilbert space $H$. Denote 
  by $H_m\subset H$ the space of $Q_m$-fixed vectors.
  Then $\cup H_m$ is dense in $H$.
  \end{proposition}

To apply this statement, we consider 
 a sequence of finite subtrees 
 $$\tochka=J_0\subset J_1\subset J_2\subset \dots\subset \T,\qquad
 \text{such that}\quad
 \cup J_m=\T$$
 and set $Q_m:=\cK(J_m)$.
 
 \sm 
 
 For any $J_1$, $J_2$ and $g\in \Hie(\T)$ we define
 the operator
 $$
 \wt\rho_{J_1,J_2}(g):H(J_2)\to H(J_1)
 $$ 
 by
 $$
 \wt\rho_{J_1,J_2}(g):= P(J_1)\rho(g)\Bigr|_{H(J_2)}.
 $$
 Clearly, 
 $$
 \wt\rho_{J_1,J_2}(h_1gh_2)=
 \wt\rho_{J_1,J_2}(g) \qquad\text{for $h_1\in \cK(J_1)$, $h_2\in \cK(J_2)$.}
 $$
 Therefore $\wt\rho(g)$ depends only on the double coset $\frg$
 containing $g$.
 
 \begin{theorem}
 	\label{th:multi}
 	 Let $\rho$ be a unitary representation
 	of the group $\Hie(\T)$.
 	The map $\wt\rho_{J_1,J_2}$ sending
 	$\frg\in \cK(J_1)\setminus \Hie(\T)/\cK(J_2)$ to
 $\wt\rho_{J_1,J_2}(\frg)$ is a representation of the category 
 $\hier$, 	 i.~e., for any
 $$\frg_1\in \cK(J_1)\setminus \Hie(\T)/\cK(J_2),\qquad
 \frg_2\in \cK(J_2)\setminus \Hie(\T)/\cK(J_3)$$
  we have
  $$
  \wt \rho _{J_1,J_2}(\frg_1)
  \,
 \wt \rho _{J_2,J_3}(\frg_2) 
 =\wt \rho _{J_1,J_3}(\frg_1\odot\frg_2).
  $$ 
\end{theorem}

The proof occupies the remaining part of this section.

\sm 

{\bf \punct Stabilizers of subtrees.%
	\label{ss:st-vertex}}
For a vertex $v\in\vert(\T)$ denote by $\cK(v)$
the stabilizer of $v$ in $\Aut(\T)$.
Denote $\cK:=\cK(\tochka)$
the stabilizer of the initial point $\tochka$. Let us describe this group.

Let $G$ be a topological  group. Consider
the countable direct product
$G^\infty =G\times G\times \dots$
equipped with the Tikhonov topology.
The infinite symmetric group $S_\infty$ acts on 
on $G^\infty$ by permutations of factors.
The {\it wreath product}
$
S_\infty \gtrdot G$
is the semi-direct product
$
S_\infty \ltimes  G^\infty
$.

\sm 

Fix $m\in \N$.  Consider a tree $T_m$,
whose vertices are enumerated by collections
$(s_1,\dots,s_l)$, where $0\le l\le m$, $s_j\in\N$,
edges have the form
$$(s_1,\dots,s_l)\,\text{---}\,(s_1,\dots,s_l,s_{l+1}).$$
The tree $T_m$ is a neighborhood of radius $m$ of $\tochka\in \T$.
see Subsect. \ref{ss:1.1}.
An element of $\Aut(T_{m+1})$ induces an automorphism of
$T_{m}$, i.~e., we have a canonical map 
$\Aut(T_{m+1})\to \Aut(T_{m})$, the kernel is a product of countable 
number 
 of copies of $S_\infty$, copies are enumerated by
vertices of $T_{m}$ of valence 1 (i.~e., vertices of the form
$(s_1,\dots,s_m)$).
The group $\Aut(T_m)$
of automorphisms of $T_m$ is
$$
\Aut(T_m)\simeq
\underbrace{S_\infty \gtrdot\biggl(S_\infty \gtrdot \Bigl(S_\infty \gtrdot\bigl(S_\infty\gtrdot\dots \bigl)
	\Bigl) \biggl)}_{\text{$m$ times}} 
$$
and  the group $\cK$ is the inverse limit,
$$
\cK\simeq \lim\limits_{\infty\longleftarrow m}
\Aut(T_m).
$$

{\bf \punct Stabilizers of finite subtrees.%
	\label{ss:st-subtree}}
Consider a subtree $J$ and its stabilizer $\cK(J)$.
Removing edges 
of $J$ from $\T$ we get a $(\T)$-covering forest,
its components $S_v$ are enumerated by vertices
$v\in \vert(J)$. Denote by $\cK(v/\!/J)$
the stabilizer of $v$ in $\Aut(S_v)$.
Clearly, 
$$
\cK(J)=\!\!\!\!\prod_{v\in\vert(J)}\!\!\!\!
\cK(v/\!/J)\simeq 
\biggl[\lim\limits_{\infty\longleftarrow m}\underbrace{S_\infty \gtrdot\biggl(S_\infty \gtrdot \Bigl(S_\infty \gtrdot\bigl(S_\infty\gtrdot\dots \bigl)
	\Bigl) \biggl)}_{\text{$m$ times}} \biggr]^{\# \vert(J)}.
$$

{\bf \punct Proof of Theorem \ref{th:multi}.}
For $m<n$,
we have a canonical epimorphism $\Aut(T_n)\to \Aut(T_m)$.
On the other hand,
there is the following (noncanonical) embedding 
$\Aut(T_m)\to \Aut(T_n)$.
Namely, $\Aut(T_m)$ is 
 the group of transformations
of the tree $T_n$ that for each vertex
$(s_1,\dots,s_m,s_{m+1},\dots s_n)$
preserve the tail $(s_{m+1},\dots s_n)$.

So we consider the groups $\Aut(T_j)$
as embedded   one to another,
$$
 \Aut(T_1)\longrightarrow \Aut(T_2)\longrightarrow
\dots \longrightarrow\cK(v).
$$
We also have
$$
\Aut(T_1)\simeq S_\infty.
$$

\begin{lemma}
	Let $\nu$ be a unitary representation of 
	$\cK=\lim\limits_{\infty\longleftarrow n} \Aut(T_n)$ in a Hilbert space $H$.
	A vector fixed by the subgroup $\Aut(T_1)\subset \cK$
	 is fixed by the whole group $\cK$.
\end{lemma}

{\sc Proof.} Denote by $Q_m$ the kernel of a map $\cK(v)\to \Aut(T_m)$. In other words
we consider automorphisms of $\T$ that fix the neighborhood
of the origin of radius $m$. Denote by $H_m\subset H$ the subspace fixed
by $Q_m$, denote $H_0:=H$. Applying Proposition \ref{pr:ism}, we get that 
$$
H=\oplus_{m=0}^\infty (H_m\ominus H_{m+1})
$$
In $H_m\ominus H_{m+1}$ we have a representation of $\Aut(T_{m+1})$.
 
Therefore, it is sufficient to prove the similar statement
for the groups $\Aut(T_m)$. Such a group
contains a chain of subgroups 
$$S_\infty=\Aut(T_1)\subset \Aut(T_2)\subset\dots \subset \Aut(T_m).$$

Consider the group
\begin{equation}
\Aut(T_2)=S_\infty \gtrdot S_\infty=S_\infty\ltimes (S_\infty)^\infty.
\label{eq:T2}
\end{equation}
The group $S_\infty$ is a type I group and all its
unitary representations are direct sums of 
irreducible representations (Lieberman, \cite{Lie}, see also \cite{Ner-book}).
Therefore  $(S_\infty)^\infty$ satisfies the same
properties; moreover
its irreducible unitary representations have the form
$$
\rho(\sigma_1, \sigma_2,\dots)=\bigotimes \rho_j(\sigma_j)
,
$$
where $\rho_j$ are irreducible unitary representations of $S_\infty$
and all but a finite number representations $\rho_j$
are trivial.
We also can write such tensor products in the form
\begin{equation}
\tau_1(\sigma_{j_1})\otimes \tau_2(\sigma_{j_2})\otimes\dots \otimes \tau_N(\sigma_{j+N})
\label{eq:tautau}
\end{equation}
omitting trivial factors and rename $\rho$ by $\tau$.
The trivial one-dimensional representation of $(S_\infty)^\infty$ 
corresponds to the empty product.  

 Consider a unitary representation of the semidirect
product (\ref{eq:T2}). Its restriction  to $(S_\infty)^\infty$
is a direct sum of irreducible representations. If we have a summand
$\bigotimes \rho_j(\sigma_j)$, then we have also
all possible (pairwise distinct) summands
$$
\tau_1(\sigma_{i_1})\otimes \tau_2(\sigma_{i_2})\otimes\dots \otimes \tau_N(\sigma_{i+N}).
$$
If the product is not empty, then we have a countable number of such summands.
An $\Aut(T_1)\simeq S_\infty$-fixed vector has components in each summand with the same norm. Therefore such components must be 0.
Thus an $\Aut(T_1)$-fixed vector 
is also $(S_\infty)^\infty$-fixed. Hence it is
$\Aut(T_2)$-fixed.
 
We apply the same argument to the group 
$$\Aut(T_3)\simeq \Aut(T_2)\ltimes (S_\infty)^{\N^2}
$$
 and its normal subgroup
$(S_\infty)^{\N^2}$. Therefore vectors fixed by $\Aut(T_2)$
are fixed  by $\Aut(T_3)$, etc.
\hfill $\square$

\sm 

Consider  a sequence $h_j\in S_\infty$ of permutations satisfying the following 
property: for each $a\in \N$ 
 the sequence $h_j$
sends $a$ to a sequence converging to $\infty$.
Then we say that $h_j$ {\it tends to infinity}.

\begin{proposition}
	If $h_j\in S_\infty$ tends to infinity, then for any unitary
	representation $\rho$ of the group $S_\infty$ the sequence
	$\rho(h_j)$ converges weakly to the projection operator
	to?? the subspace of $\cK$-fixed vectors.  
\end{proposition}

See \cite{Ner-book}, Theorem VIII.1.4.

\begin{corollary}
Let $h_j\in \Aut(T_1)\simeq S_\infty$ tends to infinity. Then for any unitary representation
$\rho$ of the group $\cK$ the sequence $\rho(h_j)$ converges to
the operator of orthogonal projection to the space of $S_\infty$-fixed vectors. 
\end{corollary}

Next, consider the group
$$
\cK(J)=\prod_{v\in\vert(J)} \cK(v/\!/J)\simeq \cK^{\# J},
$$
Consider the diagonal embedding $\upsilon_J:\Aut(T_1)\to \cK(J)^{\# J}$.

\begin{lemma}
	\label{l:d}
	Consider the diagonal embedding $d: S_\infty\to (S_\infty)^N$.
	For any unitary representation $\rho$ of group
	$(S_\infty)^N$ the subspace of $(S_\infty)^N$-fixed vectors
	coincides with the subspace of $d(S_\infty)$-fixed vectors.
\end{lemma}

\begin{corollary}
	\label{cor:last}
	Let $h_j\in \Aut(T_1)$ tend to infinity.
	Then for any unitary representation $\rho$ of $\cK(J)$ the sequence
	$\rho(\upsilon_{J}(h_j))$ converges to the projector to
	 the subspace of $\Aut(T_1)$-fixed vectors.
\end{corollary}

{\sc Proof of Lemma \ref{l:d}.}
Irreducible unitary representations of  $(S_\infty)^N$ have type $I$,   
 any irreducible representation is a tensor product
 $\rho_1(g_1)\otimes \dots \otimes \rho(g_N)$.
 Any nontrivial irreducible representation of $S_\infty$ 
 is infinite-dimensional.
 It remains to notice that
the decomposition of a tensor product $\rho_1\otimes \rho_2$
 of two nontrivial irreducible representations
of $S_\infty$ can not contain the trivial representation
(otherwise we have a Hilbert--Schmidt intertwining operator, say $A$,
from
$\rho_1$ to the representation dual to $\rho_2$;
 eigenspaces of $A^*A$ give us finite-dimensional subrepresentations
of $\rho_1$).
\hfill $\square$  

\sm

{\sc Proof of Theorem \ref{th:multi}.}
We take two double cosets $\frg_1$, $\frg_2$,
their representatives $g_1$, $g_2$, the corresponding
bi-trees $\Delta$, $\Gamma$.
Choose a sequence $h_k\in \Aut(T_1)$
tending to infinite
and take the diagonal embedding $\upsilon_{J_2}:\Aut(T_1)\to \cK(J)$
as above.
Consider the product
$
g_1\upsilon_{J_2}(h_k) g_2
$. 
Consider subtrees
$^*\!\Delta^\black_\blue$, $h\bigl(^*\Gamma^\black_\red\bigr)$.
Their intersection contains $J_2$. For sufficiently large
$k$ the transformation
$\upsilon_{J_2}(h_k)$ moves remaining pieces
of $\Delta^\black_\blue$  apart
from pieces of $\Gamma^\black_\red$,
i.~e.,
$$
^*\!\Delta^\black_\blue\bigcap\bigl(^*\Gamma^\black_\red\bigr)=J_2.
$$
By Lemma \ref{l:independent} we get gluing by $J_2$.

Next, 
\begin{multline*}
\wt\rho_{J_1,J_2} (\frg_1) \,\wt \rho_{J_2,J_3}(\frg_2)=
P(J_1)\rho(g_1) P(J_2) \cdot P(J_2)\rho(g_2) P(J_3)=
\\=
\lim_{j\to\infty} P(J_1)\rho(g_1)\, \rho(\upsilon_{J}(h_j))\, \rho(g_2) P(J_3)
= \lim_{j\to\infty} P(J_1)\rho(g_1\cdot \upsilon_{J}(h_j)\cdot g_2) P(J_3),
\end{multline*}
where $\lim_{j\to\infty}$ denotes the weak operator limit.
For sufficiently large $j$ the double coset
$\cK(J_1)\cdot g_1\cdot \upsilon_{J}(h_j)\cdot g_2\cdot\cK(J_3)$
is eventually constant and coincides with $\frg_1\odot\frg_2$.
\hfill $\square$

\section{Sphericity%
	\label{s:sphericity}}

\COUNTERS 

{\bf \punct Sphericity.}
Let $G$ be a topological group, $K$ a subgroup.
An irreducible unitary representation $\rho$ of 
 $G$ in a Hilbert space $V$ is {\it $K$-spherical}
 if the space $V^K$ of $K$-fixed vectors is one-dimensional.
 A unit vector $w\in V^K$ is called
 a {\it spherical vector}, the corresponding
 {\it spherical function}  is defined by the formula
 $$
 \Phi(g):=\la \rho(g)w,w\ra,\qquad \text{where $g\in G$.}
 $$
 The pair {\it $G\supset K$ is spherical} if for any irreducible
 unitary representation of $G$ we have $\dim V^K\le 1$.
 
 \sm 
 
 {\sc Examples.} 1) Consider the pair 
 \begin{equation}
 (S_{2\infty},S_\infty\times S_\infty)=
 S(\N\bigl|\N)
 \label{eq:NN}
 \end{equation}
 discussed in Subsect. \ref{ss:symmetric}--\ref{ss:topology}.  Consider a two-dimensional Euclidean  space
 $Y$ and two unit vectors $\xi$, $\eta\in Y$.
 Consider the countable tensor product
 $$
 \Bigl((Y,\xi)\otimes(Y,\xi)\otimes \dots \Bigr) 
 \bigotimes
 \Bigl((Y,\eta)\otimes(Y,\eta)\otimes \dots \Bigr),
 $$
  recall that the definition of a countable tensor product
  of Hilbert spaces requires a fixing of a unit vector in each 
  factor, see, e.~g., \cite{Gui}, Appendix A. Two groups $S_\infty$
  act permuting factors in big brackets,
  finitely supported permutations of $\N\sqcup\N$
  act permuting factors between brackets. Thus
  we get a representation of the group
  (\ref{eq:NN}). The vector   
 $\xi^{\otimes \infty}\otimes \eta^{\otimes\infty}$
 is spherical, the corresponding spherical function
 is
 $
 \Psi_\nu(\sigma):=
 \nu^{ m(\sigma)}
 $, where $\nu=|\la \xi,\eta\ra|^2$ and $m(\sigma)$ is the number of elements of the first 
 copy of $\N$ sent by $\sigma$ to the second copy.
 By \cite{Olsh-symm}, this one-parametric family
 of representations exhaust 
 all $S_\infty\times S_\infty$-spherical
 representations of the group (\ref{eq:NN}).
 
 \sm 
 
2) Restricting these representations to the group of spheromorphisms
 we get a one-parametric family of $\Aut(\T)$-spherical representations
 of $\Hie(\T)$. Spherical functions are given by
 the formula
 $\Phi_\nu(g)=\nu^{k(g)-1}$, where $k(g)$
 is  the number of elements in the perfect $(\T)$-covering
 forest for $g$.
\hfill $\boxtimes$

 \begin{theorem}
 	\label{th:sphericity}
 	The subgroup $\Aut(\T)$ is spherical in $\Hie(\T)$.
 \end{theorem}

For a proof we need  Olshanski's classification \cite{Ols-infinite} of unitary representation of $\Aut(\T)$.

\sm

{\bf\punct  Analog
of the complementary series for the group $\Aut(\T)$.%
\label{ss:complementary}}
The following construction arises
 to Ismagilov \cite{Ism}.
Denote by $d(v,w)$ the natural distance on the set $\vert(\T)$.
Fix real $\lambda\in [-1,1]$. Then 
$$
K_\lambda(v,w)=\lambda^{d(v,w)}
$$
is a positive definite kernel%
\footnote{See. e.g., \cite{Ner-gauss}, Section 7.1.}  on $\vert(\T)$.
Consider the Hilbert space $H_\lambda$ determined by this kernel.
In other words, we consider a Hilbert 
space $H_\lambda$ and a system of vectors $\phi_v\in H_\lambda$, where $v$
ranges in $\vert(\T)$, such that:

\sm

$\bullet$ $\la\phi_v,\phi_w\ra=\lambda^{d(v,w)}$;

\sm

$\bullet$ linear combinations of $\phi_v$ are dense in $H_\lambda$.

\sm 

For a simple explanation of the existence of this space, see, e. g., \cite{Ner3}.

\sm

A unitary representation $\Pi_\lambda$ of $\Aut(\T)$ in $H_\lambda$ 
is determined by the formula
$$\Pi_\lambda(\phi_v)=\phi_{gv}.$$

\sm

For $\lambda=0$ vectors $\phi_v$ are pairwise orthogonal
and we get the representation in $\ell^2\bigl(\vert(\T)\bigr)$.

\sm

In two cases we get degenerate constructions:

\sm

--- For $\lambda=1$
all $\phi_v$ coincide and we get the trivial one-dimensional representation
of $\Aut(\T)$.

\sm 

--- For $\lambda=-1$ we have $\phi_w=(-1)^{d(v,w)}\phi_v$
and we get a one-dimensional representation of $\Aut(\T)$.
In fact we get a homomorphism $\Aut(\T)\to\Z_2$
defined by
$$
\sigma(g)=(-1)^{d(gv,v)},
$$
where $v\in \vert(\T)$, the result does not depend on a choice of a vertex $v$.

\sm 

In nondegenerate cases finite collections of vectors $\phi_v$ are linear independent.  

\sm

{\bf\punct  Unitary  representations of $\Hie(\T)$.}

\sm 

{\sc  Cuspidal representations of $\Aut(\T)$.}
Let $J\subset \T$ be a finite subtree with $\ge 2$ vertices.  
  Denote by  $\cK(J)\subset \Aut(\T)$ the point-wise stabilizer of $J$
and by $\wt \cK(J)\subset \Aut(\T)$ the subgroup
consisting of transformations sending $J$ to itself.
Clearly, $\cK(J)$ is a normal subgroup in $\wt \cK(J)$ of finite index,
$$
\wt \cK(J)/\cK(J)\simeq \Aut(J).
$$
A {\it cuspidal representation}
of $\Aut(\T)$ is a representation induced from an
irreducible representation of a subgroup $\wt \cK(J)$
trivial on $\cK(J)$.

Notice that a cuspidal representation induced from
$\wt \cK(J)$
is a subrepresentation in $\ell^2\bigl(\Aut(\T)/\cK(J)\bigr)$.

\sm 

{\sc Classification of unitary representations.} The group
$\Aut(\T)$ has type I. Any unitary representation
of $\Aut(\T)$ is a direct integral of irreducible representations.
Any irreducible unitary representation of $\Aut(\T)$
has the form $\Pi_\lambda$ or is cuspidal.

\sm

{\bf\punct A  lemma.} 
Proof of Theorem \ref{th:sphericity} is almost identical to the proof of sphericity
for groups of spheromorphisms of Bruhat--Tits trees  in
\cite{Ner-19}. There is one place of a proof that requires
separate considerations.
 Corollary 2.5 in \cite{Ner-19} 
is based on Lemma 2.4 that
makes no sense in our case. So the corresponding statement, i.~e.,
the following Lemma \ref{l:weak0}, must be reproved.

\sm 


Let us color vertices of the tree $\T$
into two colors, say black and white,
such that each edge has ends of different colors.
Denote by $\Aut_+(\T)$ the subgroup of $\Aut(\T)$
consisting of transformations preserving the coloring. Clearly, $\Aut_+(\T)$ is a normal subgroup
in $\Aut(\T)$ of index 2.

\begin{lemma}
	\label{l:weak0}
	Let $\dots$, $a_{-1}$, $a_0$, $a_1$, $\dots$ be a two-side way in $\T$.
	Consider $h\in \Aut_+(\T)$ sending each $a_j$  to $a_{j+1}$.
	Then for any unitary representation $\rho$ of $\Hie(\T)$
	the sequence
	$\rho\bigl(h^{m}\bigr)$ weakly converges to the  operator
	of orthogonal projection
	to the space of $\Aut_+(\T)$-fixed vectors.
\end{lemma} 

{\sc Remark.} Recall a criterion of weak operator convergence.
Let $\xi_\alpha$ be a subset, whose linear combinations
are dense in a Hilbert space $H$.
Let $A_n$ be a sequence of linear operators and
the sequence $\|A_n\|$ be bounded. Then $A_n$ weakly converges
to $A$ if and only if for each $\xi_\alpha$, $\xi_\beta$
the sequence
$\la A_n \xi_\alpha, \xi_\beta\ra$
converges to $\la A \xi_\alpha, \xi_\beta\ra$.
\hfill $\boxtimes$

\sm

{\sc Proof.} It is sufficient to verify the statement for
irreducible representation of $\Aut(\T)$.

For representations $\Pi_\lambda$, where $-1<\lambda<1$,
we have
$$\la \rho(h^{m})\phi_v,\phi_w\ra\sim \lambda^{m+O(m)},
\qquad m\to+\infty
$$ 

Since any cuspidal representation is a subrepresentation
in  $\ell^2$ on some homogeneous space $\Aut(\T)/\cK(J)$, it is sufficient
to examine such representations. Denote by $e_\nu$ the standard
basis in this $\ell_2$, it is enumerated by injective maps
$J\to\T$. Clearly,
$\la \rho(h^{m}) e_\nu, e_\mu\ra$ for fixed $\mu$, $\nu$
can be nonzero only for one value of $m$.

Thus for any irreducible  infinite-dimensional representation 
of $\Aut(\T)$ the sequence $\rho(h^{m})$
weakly converges to 0. For one-dimensional representations
$\Pi_{\pm1}$ the sequence consists of unit operators.
\hfill $\square$

\sm

\section{Embeddings of $\Hie(\T)$ to infinite dimensional
	group $\GL$%
	\label{s:GL}}

\COUNTERS

{\bf \punct The infinite-dimensional group $\mathbf{GL}$.}
Let $H$ be a real infinite-di\-men\-sio\-nal Hilbert space.
Denote by $\bfO(H)$ the group of orthogonal operators
(real unitary operators) in $H$, we equip it with the weak
operator topology.
Denote by $\mathbf{GL}(H)$ the group
of all operators $g$ in a real Hilbert space $H$ admitting
a representation in the form
$g=U(1+T)$, where $U\in \bfO(H)$
and $T$ is a Hilbert--Schmidt operator%
\footnote{I. e. $\tr (T^*T)<\infty$.
	The function $(T,S)\to \tr(TS)$
	determines an inner product and a structure
	of Hilbert space on the set of all Hilbert--Schmidt
	operators. In particular this determines a topology
	on this space.
	For details, see, e. g., \cite{RS}.}. The polar decomposition
of such $g$ has the form 
$g=V \exp{S},$ where
$V\in \bfO(H)$  and $S$ is a self-adjoint Hilbert--Schmidt operator.
So we represent the space $\mathbf{GL}$
as a direct product of the group of orthogonal operators
and the space of self-adjoint Hilbert--Schmidt operators. 
This determines the
{\it Shale topology}
on $\mathbf{GL}$.


Consider the Gaussian measure $\mu$ on an extension $\Omega$ of the space
$H$ with the characteristic function%
\footnote{See, e. g., \cite{ShF} or \cite{Bog}.
	In the present paper, we do not  need precise description of this object.}
$e^{-\|h\|^2/2}$. The group
$\mathbf{GL}$ acts on $\Omega$ by linear transformations leaving
$\mu$ quasiinvariant, the orthogonal group $\bfO(H)$
preserves $\mu$.
Therefore we get a unitary representation
of $\mathbf{GL}$ in $L^2(\Omega,\mu)$, the constant function
is $\bfO(H)$-spherical. 

The group $\GLO$ is one of $(G,K)$-pairs considered
in Olshanski's theory of representations of infinite-dimensional classical
groups, see \cite{Olsh-GB}, \cite{Pick}, \cite{Ner-book}.

\sm

{\bf\punct Embeddings of $\Hie(\T)$
	to the infinite-dimensional group  $\mathbf{GL}$.}
Let $H_\lambda$ be as in Subsect. \ref{ss:complementary}.

\begin{theorem}
	For any $g\in \Hie(\T)$ there is a bounded operator
	in $H_\lambda$ such that $\sigma(g)e_v=e_{gv}$.
	Moreover, $\sigma(g)$ can be represented in the form
	$$
	\sigma(g)= U (1+Q),
	$$
	where $U$ is a unitary operator and $Q$ has finite rank.
\end{theorem}

This statement
is contained in \cite{Ner3}  for a smaller group $\Hie^\circ(\T)$,
see above Subsect. \ref{ss:frame}, 
formally we have to repeat the argumentation.

\sm 

{\sc Proof.}
It is sufficient to show that operators
$
\sigma(g)^*\sigma(g)-1
$
have finite rank. Let us evaluate the sesquilinear
form
$$
R_g(h_1,h_2)=\bigl\la (\sigma(g)^*\sigma(g)-1)h_1,h_2\bigr\ra
=\la \sigma(g) h_1, \sigma(g) h_2\ra- \la h_1,h_2\ra.
$$
For a subtree $S\subset \T$ denote by 
$H_S$ the subspace in $H$ generated by $\phi_v$, where
$v\in \vert(S)$. Consider the perfect $(\T)$-forest
 $S_1$, \dots, $S_N$ for $g$. 
For each pair $S_\alpha$, $S_\beta$ we take 
nearest vertices $u_{\alpha\beta}\in S_\alpha$, $v_{\alpha\beta}\in S_\beta$.
Also we take nearest vertices $w_{\alpha\beta}\in gS_\alpha$,
$z_{\alpha\beta}\in gS_\beta$.

Clearly, the form $R_g$ is zero on $H_{S_\alpha}\times H_{S_\alpha}$
for all $S_\alpha$. 

Let $\beta\ne\alpha$.
Consider separately forms $\la h_1,h_2\ra$
and $\la \sigma(g) h_1,\sigma(g)h_2\ra$.
Let $p\in 
\vert(S_\alpha)$, $q\in \vert(S_\beta)$. Then
\begin{multline*}
\la \phi_p, \phi_q\ra=
\lambda^{d(p,q)}=\lambda^{d(p, u_{\alpha\beta})
	+d(u_{\alpha\beta},v_{\alpha\beta})+d(v_{\alpha\beta},q)
}
=\\=
\la \phi_p,\phi_{u_{\alpha\beta}}\ra
\cdot \lambda^{d(u_{\alpha\beta},v_{\alpha\beta})}
\cdot
\la \phi_{v_{\alpha\beta}}, \phi_q\ra.
\end{multline*}
Therefore, for $h_1\in H_\alpha$, $h_2\in H_\beta$ we have
$$
\la h_1,h_2\ra = 
\lambda^{d(u_{\alpha\beta},v_{\alpha\beta})}
\la h_1,\phi_{u_{\alpha\beta}}\ra \la \phi_{v_{\alpha\beta}},h_2\ra.
$$
Thus the form $\la h_1,h_2\ra$ has rank 1 on $H_{S_\alpha}\times H_{S_\beta}$.

The same argument shows that
$$
\la \sigma(g) h_1,\sigma(g) h_2\ra = 
\lambda^{d(w_{\alpha\beta},z_{\alpha\beta})}
\la h_1,\phi_{w_{\alpha\beta}}\ra \la \phi_{z_{\alpha\beta}},h_2\ra.
$$
Therefore the form $\la \sigma(g) h_1,\sigma(g) h_2\ra$
also has rank 1 on $H_{S_\alpha}\times H_{S_\beta}$.
\hfill $\square$

 \tt
\noindent
Yury Neretin\\
Wolfgang  Pauli Institute/c.o. Math. Dept., University of Vienna \\
\&Institute for Theoretical and Experimental Physics (Moscow); \\
\&MechMath Dept., Moscow State University;\\
\&Institute for Information Transmission Problems;\\
yurii.neretin@math.univie.ac.at
\\
URL: http://mat.univie.ac.at/$\sim$neretin/

\begin{thebibliography}{cc}
	
	\bibitem{Bog}
Bogachev, V. I.
{\it 
Gaussian measures.} Providence,  American Mathematical Society (AMS), 1998.
	
	
	\bibitem{Bur}	
	Burillo J., Cleary S., Stein M., Taback J. {\it Combinatorial and metric properties of Thompson’s
		group $T$}. Trans. Amer. Math. Soc. 361 (2009), no. 2, 631-652.
	
		\bibitem{Fos}
	Fossas A.
	{\it $\PSL(2,\Z)$ as a non-distorted subgroup of Thompson's group $T$.}
	Indiana Univ. Math. J. 60 (2011), no. 6, 1905-1925. 
	
	 \bibitem{GS}
	Ghys \'E., Sergiescu V.
	{\it Sur un groupe remarquable de diff\'eomorphismes du cercle.} 
	Comment. Math. Helv. 62 (1987), no. 2, 185-239.
	
	\bibitem{Gui}
	Guichardet A.
	{\it Symmetric Hilbert spaces and related topics. 
		Infinitely divisible positive definite functions. Continuous products and tensor products. Gaussian and Poissonian stochastic processes.}
	Lect. Not.  Math., Vol. 261. Springer-Verlag, Berlin-New York, 1972.
	
	\bibitem{Imb}
	Imbert, M.
	{\it Sur l'isomorphisme du groupe de Richard Thompson avec le groupe de Ptol\'em\'ee.} 
	In {\it Geometric Galois actions}, V. 2 (eds. L. Schneps, P. Lochak.), 313-324, 
	Cambridge Univ. Press, Cambridge, 1997. 
	
	
	

	
\bibitem{Ism}
Ismagilov, R. S.,
{\it Elementary spherical functions on the groups $\SL(2,P)$ over
a field $P$, which is not locally compact with respect to
the subgroup of matrices with integral elements.} 
Mathematics of the USSR-Izvestiya, 1967, 1:2, 349-380 	
	
\bibitem{Kech}
Kechris, A. S.
{\it Classical descriptive set theory}. 
 Berlin: Springer-Verlag. xx, 402 p. (1995).
 
 \bibitem{KR}		
 Kechris 	A. S.,  Rosendal C., {\it Turbulence, amalgamation, and generic automorphisms of homogeneous structures}, Proc. Lond. Math. Soc. (3), 94:2 (2007), 302--350.	
	
 \bibitem{Lie}
Lieberman  A., {\it The structure of certain unitary representations of infinite
	symmetric groups}, Trans. Amer. Math. Soc., 164 (1972), 189-198.	
	
\bibitem{Lusin0}	
Lusin, N. {\it Sur un exemple arithm\'etique d'une fonction ne faisant pas partie de la
classification de M. Ren\'e Baire}.  Compt. Rend. 182, 1521-1522 (1926).

	
\bibitem{Lusin}
Lusin, N.
{\it Le\c{c}ons sur les ensembles analytiques et leurs applications.}
 Paris, Gauthier-Villars (1930).

	  \bibitem{Ner1}
Neretin  Yu.A.  {\it Unitary representations of the groups of diffeomorphisms of the $p$-adic projective line,}
Functional Anal. Appl. 18 (1984) 345-346). 

\bibitem{Ner2}
Neretin Yu.A. {\it Combinatorial analogues of the group of diffeomorphisms of the circle,}
Russian Acad. Sci. Izvestiya. Math. 41 (2) (1993) 337-349.


\bibitem{Ner3}
Neretin Yu. A. {\it Groups of hierarchomorphisms of trees and related Hilbert spaces.} J. Funct. Anal. 200 (2003), no. 2, 505-535. 

 \bibitem{Ner-book}
Neretin Yu. A. {\it Categories of symmetries and infinite-dimensional groups.}  The Clarendon Press, Oxford University Press, New York, 1996.

\bibitem{Ner-gauss}	 
Neretin Yu. A., {\it Lectures on Gaussian integral operators and classical groups}, EMS Ser. Lect. Math., Eur. Math. Soc. (EMS), Z\"urich, 2011.

 \bibitem{Ner-symm}
Neretin Yu. A. {\it Infinite symmetric groups and combinatorial constructions of topological field theory type.}
Russian Math. Surveys 70 (2015), no. 4, 715-773.

\bibitem{Ner-19}
Neretin Yu. A.
{\it On spherical unitary representations of groups of spheromorphisms of Bruhat--Tits trees.} Preprint, arXiv:1906.12197 

 \bibitem{Olsh-trees}
Olshanski  G.I.
{\it Classification of the irreducible representations of the automorphism groups of Bruhat-Tits trees},
Funct. Anal. Appl. 11(1) (1977), 26-34.	

\bibitem{Ols-infinite}
 Olshanskii  G. I., {\it New “large” groups of type I},
  J. Soviet Math., 18:1 (1982), 22-39.
  
   \bibitem{Olsh-GB}
  Olshanski G. I.
  {\it Unitary representations of infinite-dimensional pairs $(G,K)$ and the formalism of R. Howe.}
  In Zhelobenko D. P., Vershik A. M. (eds.) {\it Representation of Lie groups and related topics}, 269-463, 
  Adv. Stud. Contemp. Math., 7, Gordon and Breach, New York, 1990.
  
   \bibitem{Olsh-symm}
  Olshanski G. I. {\it Unitary representations of $(G,K)$-pairs that are connected with the infinite symmetric group $S(\infty)$}.
  Leningrad Math. J. 1 (1990), no. 4, 983-1014.
  
  \bibitem{Pick}
  Pickrell, D. {\it  Separable representations for automorphism 
  	groups of infinite symmetric spaces.}  J. Funct. Anal.  90  (1990),  no. 1, 1-26.

\bibitem{RS}
Reed, M.,
 Simon, B.
{\it Methods of modern mathematical physics. I: Functional analysis.}
New York-London: Academic Press, 1972.

\bibitem{ShF}
 Shilov, G. E.; Fan Dyk Tin {\it Integral, measure and derivative on linear
 space}. Nauka, Moscow 1966.
(Russian)
	
\end{thebibliography}
\end{document}